\documentclass[a4paper,12pt,leqno]{amsart}
\usepackage{amsmath}
\usepackage{cases}
\usepackage{bm}
\usepackage{amssymb}
\usepackage{indentfirst}
\usepackage{mathrsfs}
\usepackage{setspace}
\usepackage{geometry}
\usepackage{tikz}
\usepackage{tikz}
\RequirePackage{amsmath,amsthm,amssymb}
\theoremstyle{plain}
\providecommand{\theoremname}{Theorem}%

\providecommand{\axiomname}{Axiom}%

\providecommand{\lemmaname}{Lemma}%

\providecommand{\corollaryname}{Corollary}%

\providecommand{\assertionname}{Assertion}%

\providecommand{\propositionname}{Proposition}%

\providecommand{\conjecturename}{Conjecture}%

\theoremstyle{definition}
\providecommand{\definitionname}{Definition}%

\providecommand{\examplename}{Example}%

\theoremstyle{remark}
\providecommand{\remarkname}{Remark}%

\begin{document}
\begin{spacing}{1.5}

 \def\Cal{\mathcal}
\def\bold{\mathbf}
\def\ca{\mathcal A}
\def\cdz{\mathcal D_0}
\def\cd{\mathcal D}
\def\cdo{\mathcal D_1}
\def\bold{\mathbf}
\def\l{\lambda}
\def\le{\leq}

\def\ll{\underset {L}{\leq}}
\def\rl{\underset {R}{\leq}}
\def\lr{\rl}
\def\lrl{\underset {LR}{\leq}}
\def\llr{\lrl}
\def\el{\underset {L}{\sim}}
\def\er{\underset {R}{\sim}}
\def\elr{\underset {LR}{\sim}}
\def\ds{\displaystyle\sum}

\title[Based Ring of Two-sided Cells ]
{The Based Rings of Two-sided cells in an Affine Weyl group of type $\tilde B_3$, IV}

\author[Y. Qiu]{Yannan Qiu$^{*}$}

\begin{abstract}
We compute the based ring of two-sided cell corresponding to the unipotent class
 in $Sp_6(\mathbb C)$ with Jordan blocks (21111). The results also verify Lusztig's conjecture on the structure of the based rings of the two-sided cells of an affine Weyl group.
\end{abstract}

\maketitle

Let $G$ be $Sp_6(\mathbb C)$ and $W$ be the extended affine Weyl group attached to $G$.  We are concerned with the based rings of two-sided cells of   $W$. In   previous papers [QX2, QX1,Q] we computed the  based rings  of the two-sided cells of $W$  with values 2, 3, 4 of Lusztig's $a$-function.  In this paper we compute the based ring of two-sided cell $c$ of $W$ with value 6 of Lusztig's $a$-function. The result also verifies Lusztig's conjecture on the structure of the based ring of a two-sided cell for $c$. For this two-sided cell, the validity of Lusztig's conjecture on the based rings is already included in main theorem in [BO]. Here we construct the bijection in Lusztig's conjecture explicitly so that the result in this paper can be used for  computing irreducible representations of affine Hecke algebras of type $\tilde B_3$.

The contents of the paper are as follows. Section 1 is devoted to set up and notations. In section 2 we  recall some results on cells of $W$, which are due to J. Du, and establish some conclusions for  the two-sided cell $c$. Sections 3 is devoted to computing the based ring  of the two-sided cell $c$. Under a bijection of Lusztig, this two-side cell corresponds to the unipotent class  in $Sp_6(\mathbb C)$ with Jordan blocks (21111).

\def\Cal{\mathcal}
\def\bold{\mathbf}
\def\ca{\mathcal A}
\def\cdz{\mathcal D_0}
\def\cd{\mathcal D}
\def\cdo{\mathcal D_1}
\def\bold{\mathbf}
\def\l{\lambda}
\def\le{\leq}

\def\ll{\underset {L}{\leq}}
\def\rl{\underset {R}{\leq}}
\def\lr{\rl}
\def\lrl{\underset {LR}{\leq}}
\def\llr{\lrl}
\def\el{\underset {L}{\sim}}
\def\er{\underset {R}{\sim}}
\def\elr{\underset {LR}{\sim}}
\def\ds{\displaystyle\sum}

\section{Set up and notations}

In this section we fix some notations and recall some basic facts. We refer to [KL, L2, L3, L4, QX2] for more details.

\medskip

 \noindent{\bf 1.1. Extended affine Weyl groups and their Hecke algebras} \
 Let $G$ be a simply connected, almost simple complex algebraic group and $T$ a maximal torus of $G$. Let $\Phi\subset X=Hom(T, \mathbb C^*)$ be the root system, $P\subset X$ the root lattice. Let $\Delta$ be the set of simple roots in $\Phi$ and $\Phi^+$ be the set of positive roots in $\Phi$.The Weyl group $W_0=N_G(T)/T$ of $G$ acts on $X$ in a natural way and this action is stable on $P$ and $\Phi$. Thus we can form the affine Weyl group $W'=W_0\ltimes P$ which is a normal subgroup of the extended affine Weyl group $W=W_0\ltimes X$. There exists a finite abelian subgroup $\Omega$ of $W$ such that $W=\Omega\ltimes W'$. The set of simple reflections of $W'$ is denoted by $S$. We shall denote the length function of $W$ by $l$ and use $\leq$ for the
Bruhat order on $W$. We also often write $y<w$ or $w>y$ if $y\le w$ and $y\ne w$.

Let $H$ be the Hecke algebra  of $(W,S)$ over $\mathcal A=\mathbb
C[q^{\frac 12},q^{-\frac 12}]$ with
parameter $q$. Let  $\{T_w\}_{w\in W}$ be its standard basis. Then we have $(T_r-q)(T_r+1)=0$ for $r\in S$ and $T_wT_u=T_{wu}$ if $l(wu)=l(w)+l(u)$.
Let
$C_w=q^{-\frac {l(w)}2}\sum_{y\le w}P_{y,w}T_y,\ w\in W$ be the
Kazhdan-Lusztig basis of $H$, where $P_{y,w}$ are the
Kazhdan-Lusztig polynomials.  The degree of $P_{y,w}$ is less than
or equal to $\frac12(l(w)-l(y)-1)$ if $y< w$ and  $P_{w,w}=1$. Convention: set $P_{y,w}=0$ if $y\not\le w$.

\def\vp{\varphi}
\def\st{\stackrel}
\def\sc{\scriptstyle}

 We should remark here that our notation $C_w$ stands for $C'_w$ in [KL]. The reason we use the elements $C'_w$ in [KL] since the positivity in multiplications of the elements is simpler in writing and in practical computation.

The following  formulas for computing $C_w$ (see [L2]) will be used in Section 3. (Note the notation $C_w$ here stands for $C'_w$ in [KL].)

\medskip

We shall write $y\prec w$ if $\mu(y,w)\ne 0$.  We have

(a)  Let $y,w\in W$ and $r\in S$ be such that $y<w,\ rw<w,$ and $ ry>y$. Then $y\prec w$ if and only if $w=ry$. Moreover this implies that $\mu(y,w)=1$.

(b) Let $w\in W$ and $r\in S$. Then
$$\begin{aligned} C_rC_w=\begin{cases}\displaystyle (q^{\frac12}+q^{-\frac12})C_w,\quad &\text{if\ }rw<w,\\
\displaystyle C_{rw}+\sum_{\st  {z\prec w}{rz<z}}\mu(z,w)C_z,\quad&\text{if\ }rw\ge w.\end{cases}\end{aligned}$$
$$\begin{aligned} C_wC_r=\begin{cases}\displaystyle (q^{\frac12}+q^{-\frac12})C_w,\quad &\text{if\ }wr<w,\\
\displaystyle C_{wr}+\sum_{\st  {z\prec w}{zr<z}}\mu(z,w)C_z,\quad&\text{if\ }wr\ge w.\end{cases}\end{aligned}$$

\medskip

\noindent{\bf 1.2. Cells of affine Weyl groups}  We refer to [KL], [L2] for definition of left cells, right cells and two-sided cells of $W$.

For $h,\, h'\in H$ and $x\in W$, write
$$ hC_x =\sum_{y\in W}a_yC_y,\quad
 C_xh =\sum_{y\in W}b_yC_y,\quad
  hC_xh'=\sum_{y\in W}c_yC_y,\quad a_y,b_y, c_y\in \mathcal A.$$
  Define $y\ll x$ if $a_y\ne 0$ for some $h\in H$, $y\rl x$ if $b_y\ne 0$ for some $h\in H$, {and} $y\lrl x$ if $c_y\ne 0$ for some $h,h'\in H$.

  We write $x\el y$ if $x\ll y\ll x$,  $x\er y$ if $x\rl y\rl x$, and $x\elr y$ if $x\lrl y\lrl x$. Then $\el,\ \er,\ \elr$ are equivalence relations on $W$. The equivalence classes are called left cells, right cells, and two-sided cells of $W$ respectively. Note that if $\Gamma$ is a left cell of $W$, then $\Gamma^{-1}=\{ w^{-1}\,|\, w\in\Gamma\}$ is a right cell.

\medskip

For $w\in W$, set $R(w)=\{r\in S\,|\, wr\le w\}$ and $L(w)=\{r\in S\,|\, rw\le w\}.$ Then we have (see [KL])

\medskip

(a) $R(w)\subset R(u)$ if $u\ll w$ and $L(w)\subset L(u)$ if $u\rl w.$ In particular, $R(w)= R(u)$ if $u\el w$ and $L(w)= L(u)$ if $u\er w.$

\medskip

\noindent{\bf 1.3.   $*$-operations}\ \ The $*$-operation introduced in [KL] and   generalized in [L2] is a useful tool in the theory of cells of Coxeter groups.

Let $r,t$ be simple reflections in $S$ and assume that $rt$ has order  $m\ge 3$. Let $w\in W$ be such that $sw\ge w,\ tw\ge w$. The $m-1$ elements $rw,\ trw,\ rtrw,\ ...,$ are called a left string (with respect to $\{r,t\})$, and the $m-1$ elements $tw,\ rtw,\ trtw,\ ..., $ are also called a left string (with respect to $\{r,t\})$. Similarly we define right strings (with respect to $\{r,t\}$). Then (see [L2])

\medskip

(a) A left string in $W$ is contained in a left cell of $W$ and a right string in $W$ is contained in a right cell of $W$.

\medskip
 Assume that $x$ is in a left (resp. right ) string (with respect to $\{r,t\})$ of length $m-1$ and is the $i$th element of the left (resp. right) string,  define ${}^*x$ (resp. $x^*$) to be the $(m-i)$th element of the string, where $*=\{r,t\}$. The following result is proved in [X2].

 \medskip

 (b) Let $x$ be in $W$ such that $x$ is in a left string with respect to $*=\{r,t\}$ and is also in a right string with respect to $\star=\{r',t'\}$. Then ${}^*x$ is in a right string with respect to $\{r',t'\}$ and $x^\star$ is in a left string with respect to $\{r,t\}$. Moreover ${}^*(x^\star)=({}^*x)^\star$. We shall write ${}^*x^\star$
 for ${}^*(x^\star)=({}^*x)^\star$.

\medskip

The following result is due to Lusztig [L2].

\medskip

(c) Let $\Gamma$ be a left cell of $W$ and an element $x\in\Gamma$ is in a right string $\sigma_x$ with respect to $*=\{r,t\}$. Then any element $w\in\Gamma$ is in a right string $\sigma_w$ with respect to $*=\{r,t\}$. Moreover $\Gamma^*=\{w^*\,|\, w\in \Gamma\}$ is a left cell of $W$ and $\mathscr C=\displaystyle \left(\cup_{w\in\Gamma}\sigma_w\right)-\Gamma$ is a union of at most $m-2$ left cells.

\medskip

Following Lusztig [L1] we set $\tilde\mu(y,w)=\mu(y,w)$ if $y\leq w$ and $\tilde\mu(y,w)=\mu(w,y)$ if $w<y$. For convenience we also set $\tilde\mu(y,w)=0$ if $y\nleq w$ and $w\nleq y$. Assume that $x_1,x_2,...,x_{m-1}$  and $y_1,y_2,...,y_{m-1}$ are two left strings with respect to $*=\{r,t\}$. Define
$$a_{ij}=\begin{cases}\tilde\mu(x_i,y_j),\quad &\text{if\ } \{r,t\}\cap L(x_i)=\{r,t\}\cap L(y_j),\\
0,\quad &\text{otherwise}.\end{cases}$$
Lusztig proved the following identities (see Subsection 10.4 in [L1]).

\medskip

(d) If $m=3$, then $a_{11}=a_{22}$ and $a_{12}=a_{21}$.

\medskip

(e) If $m=4$, then
$$\begin{aligned} a_{11}=a_{33},\ a_{13}=a_{31},\ a_{22}=a_{11}+a_{13},\ a_{12}=a_{21}=a_{23}=a_{32}.\end{aligned}$$

\medskip

\noindent{\bf 1.4. Lusztig's $a$-function}\quad
For $x,y\in W$, write $$C_xC_y=\sum_{z\in
W}h_{x,y,z}C_z,\qquad h_{x,y,z}\in \mathcal A
$$
 Following Lusztig ([L2]), we define
 $$a(z)={\rm max}\{\text{deg}\ h_{x,y,z}| x,y\in W\}.$$

 Springer showed that $l(z)\ge a(z)$ (see [L3]). Let $\delta(z)$ be
the
 degree of $P_{e,z}$, where $e$ is the neutral element of $W$.
 Then actually one has $l(z)-a(z)-2\delta(z)\ge 0$ (see [L3]). Set
 $$\cd =\{z\in W\ |\ l(z)-a(z)-2\delta(z)=0\}.$$
The elements of $\cd$ are involutions, called distinguished involutions of
 $(W,S)$ (see [L3]).

 The following properties  are proved in [L1].

(a) $a(x)\ge a(y)$ if $x\lrl y$. In particular, $a(x)=a(y)$ if $x\elr y$.

(b) If $h_{x,y,z}\ne 0$, then $z\lr x$ and $z\ll y$. In particular, $a(z)>a(x)$ if $z\not \er x$, and $a(z)>a(y)$ if $z\not \el y$.

Let $i$ be a nonnegative integer. Using (a), (b) and the definitions in Subsection 1.2, we see that  the $\mathcal A$-submodule $I$ of $H$ spanned by all $C_w$ with $a(w)\ge i$ is a two-sided ideal of $H$. This fact will be used in Section 3, 4 and 5.

Following Lusztig,
 we define  $\gamma_{x,y,z}$ by the following formula,
 $$h_{x,y,z}=\gamma_{x,y,z}q^{\frac {a(z)}2}+
 {\rm\ lower\ degree\ terms}.$$
  We remark that the notation $\gamma_{x,y,z}$ here stands for Lusztig's original notation $\gamma_{x,y,z^{-1}}$ in   [L3].
The following properties are due to Lusztig [L3] except  (f) (which is trivial) and (g) (proved in [X2]).

 (c) $\gamma_{x,y,z}\ne 0\Longrightarrow x\el y^{-1},\ y\el z,\ x\er z.$

 (d) $\gamma_{x,y,z}=\gamma_{y,z^{-1},x^{-1}}=\gamma_{z^{-1}, x,y^{-1}}$.

 (e) $\gamma_{x,y,z}=\gamma_{y^{-1},x^{-1},z^{-1}}.$

 (f) If $\omega,\tau\in W$ has length 0, then
 $$\gamma_{\omega x,y\tau,\omega z\tau}=\gamma_{x,y,z},\ \
 \gamma_{ x\omega,\tau y, z}=\gamma_{x,\omega\tau y,z}.$$

 (g) Let $x,y,z\in W$ be such that (1) $x$ is in a left string with respect to $*=\{r,t\}$ and also in a right string with respect to $\#=\{r',t'\}$, (2) $y$ is in a left string with respect to $\#=\{r',t'\}$ and also in a right string with respect to $\star=\{r'',t''\}$, (3) $z$ is in a left string with respect to $*=\{r,t\}$ and also  in a right string with respect to $\star=\{r'',t''\}$. Then
 $$\gamma_{x,y,z}=\gamma_{{}^*x^\#,{}^\#y^\star,{}^*z^\star}.$$

\def\tt{\tilde T}

\medskip

 \noindent{\bf 1.5.} Assume $r, t\in S$ and $rt$ has order 4. Let $w, u, v, w', u', v'$ be in $W$ such that $l(rtrtw)=l(w)+4$, $l(rtrtv)=l(v)+4$, $l(u'rtrt)=l(u')+4$ and $l(v'rtrt)=l(v')+4$. We have (see [X2 1.6.3])
\begin{itemize}
\item[(a)] $\gamma_{trw, u, trv}= \gamma_{rw, u, rv}+ \gamma_{rw, u, rtrv},$
\item[(b)]$\gamma_{w', u't, v'trt}+\gamma_{w', u't, v't}=\gamma_{w', u'tr, v'tr},$
\end{itemize}

\medskip

\def\ll{\underset {L}{\leq}}
\def\rl{\underset {R}{\leq}}

\def\lrl{\underset {LR}{\leq}}
\def\llr{\lrl}
\def\el{\underset {L}{\sim}}
\def\er{\underset {R}{\sim}}
\def\elr{\underset {LR}{\sim}}
\def\ds{\displaystyle\sum}

\def\vp{\varphi}
\def\st{\stackrel}
\def\sc{\scriptstyle}

 \noindent{\bf 1.6. The based ring of a two-sided cell}\quad For each two-sided cell $c$ of $W$, let $J_c$ be the free $\mathbb Z$-module with a basis $t_w,\ w\in c$. Define
  $$t_xt_y=\sum_{z\in c}\gamma_{x,y,z}t_z.$$
  Then $J_c$ is an associative ring with unit $\sum_{d\in\cd\cap c}t_d.$
 Then $J_c$ is an associative ring with unit $\sum_{d\in\cd\cap c}t_d.$   The ring $J=\bigoplus_{c}J_\mathfrak c$ is a based ring with unit $\sum_{d\in\cd}t_d$. {We have $J_c=span\{t_w\ |\ w\in c\}=(\sum_{d\in\cd\cap c}t_d)J(\sum_{d\in\cd\cap c}t_d)$. If  $\Gamma\subset c$ is a left cell, and $d_\Gamma$ is the unique distinguished involution in $\Gamma$, then
$$J_{\Gamma\cap\Gamma^{-1}}=t_{d_\Gamma}Jt_{d_\Gamma}\subset J_c$$ is also a based ring, with identity element $t_{d_\Gamma}$.
If $\Gamma, \Theta$ are two left cells, then $t_{d_\Gamma}Jt_{d_{\Theta}}=J_{\Gamma\cap{\Theta}^{-1}}$ is a based $t_{d_\Theta}Jt_{d_\Theta}-t_{d_\Gamma}Jt_{d_\Gamma}$-bimodule.}
Sometimes $J$ is called the  asymptotic Hecke algebra since Lusztig established an injective $\ca$-algebra homomorphism (see [L3])
  $$  \phi: H \to J\otimes\ca,\quad
  C_x \mapsto\sum_{\st {\st {d\in\cd}{w\in W}}{w\el d}}h_{x,d,w}t_w.$$

  \medskip

 \noindent{\bf 1.7.} For any weight $x\in X$, let $y$ and $z$ be dominate weights such that $x=yz^{-1}$ (the operation in $X$ is written multiplicatively). Set $\theta_x=q^{-l(y)/2}T_y(q^{-l(z)/2}T_z)^{-1}=\tilde{T_y}\tilde{T_z}^{-1}$. Then $\theta_x$ is well defined and independent of the choice of $y$ and $z$. Moreover, for any $x,x'\in X$ we have $\theta_{xx'}=\theta_x\theta_{x'}$. The following result is due to J. Bernstein:

  \medskip

  (a) For a dominant weight $x\in X$, the element
  $$S_x=\sum_{y\in X}d(y,x)\theta_y$$
  is in the center of $H$ and the  center of $H$ is spanned by all $S_x,\ x\in X$ is dominant, where $d(y,x)$ is the dimension of weight space $V(x)_y$ of an irreducible $G$-module $V(x)$ with highest weight $x$.

  \medskip

  The following result is proved in [L4].

  (b) For a dominant weight $x$, $\phi(S_x)$ is in the center of $J\otimes \mathcal A$.

  \medskip

  Let $\Gamma$ be a left cell of $W$  and $x\in X$ be dominant. Then for any $y\in \Gamma\cap\Gamma^{-1}$, since $S_xC_y=C_yS_x$, we have

  (c) $S_xC_y=\sum_{z\in\Gamma \cap\Gamma^{-1}}\zeta_{x,y,z}C_z+\heartsuit,$ where $\zeta_{x,y,z}\in\mathcal A$ and $\heartsuit$ is a linear combination of some $C_u$ with $a(u)>a(y)$.

  \medskip

  \noindent{\bf 1.8.}  Assume further that $G$ is a simply connected simple algebraic group over $\mathbb C$ and $W$ be the extended affine Weyl group attached to $G$. In [L4] Lusztig establishes a bijection between the set of two-sided cells of  $W$ and the set of unipotent classes of $G$. (In the case $Sp_6(\mathbb C)$ concerned in subsequent sections, the bijection was already  established in [D].)

  For each two-sided cell $c$ of $W$, let $u$ be a unipotent element in the unipotent class corresponding to $c$. Let $f: SL_2(\mathbb C)\to G$ be algebraic group homomorphism such that $f$ sends $\displaystyle \begin{pmatrix}1&1\\0&1\end{pmatrix}$ to $u$. Then $F_c=\{g\in G\,|\, gf(x)=f(x)g,\ \forall x\in SL_2(\mathbb C)\}$ is the  maximal reductive subgroup of the centralizer of $u$ in $G$.

  Let $\mathbf J_c=J_c\otimes_{\mathbb Z}\mathbb C$. According to [L4, 4.3], the isomorphism classes of simple $\mathbf J_c$-modules are one to one corresponding to the $F_c$-conjugate classes of pairs $(s,\rho)$, where $s\in F_c$ is semisimple, $\rho$ is a representation of $Z_{F_c}(s)/Z_{F_c}(s)^\circ$ appears in $H^*(\mathcal B^s_u,\mathbb C)$, here $\mathcal B_u^s$ is the variety of all Borel subgroup containing $s$ and $u$.

  Let $E_{s,\rho}$ be a simple $\mathbf J_c$-module corresponding to the pair $(s,\rho)$. Through the homomorphism $$\phi_c:H\to J_c\otimes_{\mathbb Z}\mathcal A=J_c\otimes_{\mathbb Z}\mathbb C[q^{\frac{1}{2}}, q^{-\frac{1}{2}}],\quad C_x \mapsto\sum_{\st {\st {d\in\cd\cap c}{w\in W}}{w\el d}}h_{x,d,w}t_w,$$  $E_{s,\rho}$ is endowed with an $H$-module structure, denoted by ${}^\phi E_{s,\rho}$.

  According to Theorem 4.2 in [L4], we have

  \medskip
  \def\tr{\text{tr}}

  (a) $\phi(S_x)$ acts on ${}^\phi E_{s,\rho}$ by scalar $\displaystyle\tr(ss_{q^{-1/2}},V(x))$, where $s_{q^{-1/2}}=\displaystyle\begin{pmatrix}q^{-\frac12}&0\\0&q^{\frac12}\end{pmatrix}$. (Note that $C_w$ in this paper is the $C'_w$ in [KL].)

\medskip

  \noindent{\bf 1.9. Lusztig's conjecture on the structure of $J_c$}
  
  {\bf Conjecture} (Lusztig [L4]): Keep the notations above. Then there exists a finite set $Y$ with an algebraic action of $F_c$ and a  bijection
  $$\pi: c\to \text{{the set of} isomorphism classes of irreducible $F_c$-vector bundles  on}\ Y\times Y.$$
  such that

  (i)  The bijection $\pi$ induces a based ring isomorphism (see [L4] for definition)
  $$\pi: J_c\to K_{F_c}(Y\times Y),\ \ t_x\mapsto \pi(x).$$

  (ii) $\pi(x^{-1})_{(a,b)}=\pi(x)_{(b,a)}^*$ is the dual representation of $\pi(x)_{(b,a)}.$

  \medskip
  
  {\bf Remark:} The conjecture was verified for various cases: $GL_n(\mathbb C),$ $SL_n(\mathbb C),$ rank 2 cases, the lowest two-sided cell, the second highest two-sided cell (see [X2, X1, X3], and a weak form (see [BO]).
  
  For the two-sided cell corresponding to the unipotent element in $Sp_6(\mathbb C)$ with three equal Jordan blocks, it is showed that the conjecture needs modification (see [QX1, BDD]). In [QX1], a conjectural geometric realization of $J_c$ for a two-sided cell of an affine Weyl group is stated. This conjectural geometric realization of $J_c$ is proved by R. Bezrukavnikov, I. Karpov and V. Krylov [BKK] and O. Propp [P]. In [P], Propp also partially resolved a  modification of Lusztig's conjecture on the structure of $J_c$ stated in [QX2].

  \section{Cells in an extended affine Weyl group of type $\tilde B_3$}

  In this section $G=Sp_6(\mathbb C)$, so that the extended affine Weyl group $W$ attached to $G$ is of type $\tilde B_3$. The left cells and two-sided cells are described by J. Du (see [D]). We recall his results.

    \noindent{\bf 2.1. The Coxeter graph of $W$}. As usual, we number the 4 simple reflections $r_0,\ r_1,\ r_2,\ r_3$ in $W$ so that
  \begin{alignat*}{2} &r_0r_1=r_1r_0,\quad r_0r_3=r_3r_0,\quad r_1r_3=r_3r_1,\\
  &(r_0r_2)^3=(r_1r_2)^3=e,\quad (r_2r_3)^4=e,\end{alignat*}
  where $e$ is the neutral element in $W$. The relations among the simple reflections can be read through the following Coxeter graph:

 
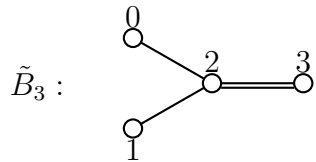
\begin{center}
  \begin{tikzpicture}[scale=.6]
    \draw (-1,0) node[anchor=east]  {$\tilde B_3:$};
    \draw[thick] (2 cm,0) circle (.2 cm) node [above] {$2$};
    \draw[xshift=2 cm,thick] (150:2) circle (.2 cm) node [above] {$0$};
    \draw[xshift=2 cm,thick] (210:2) circle (.2 cm) node [below] {$1$};
    \draw[thick] (4 cm,0) circle (.2 cm) node [above] {$3$};
    \draw[xshift=2 cm,thick] (150:0.2) -- (150:1.8);
    \draw[xshift=2 cm,thick] (210:0.2) -- (210:1.8);
    \draw[thick] (2.2,0) --+ (1.6,0);
    \draw[thick] (2.2,-0.1) --+ (1.6,0);
 \end{tikzpicture}
\end{center}

  There is a unique nontrivial element $\tau $ in $W$ with length 0. We have $\tau^2=e,\ \tau r_0\tau=r_1,\ \tau r_i\tau=r_i$ for $i=2,3.$ Note that $r_1,r_2,r_3$ generate the Weyl group $W_0$ of type $B_3$ and $r_0,r_1,r_2,r_3$ generate an affine Weyl group $W'$ of type $\tilde B_3$. And $W$ is generated by $\tau,\ r_0,r_1,r_2,r_3$.

  \noindent{\bf 2.2. Cells in $W$}\quad According to [D], the extended affine Weyl group $W$ attached to  $Sp_6(\mathbb C)$ has 8 two-sided cells:
    $$A,\quad B,\quad C,\quad D,\quad E, \quad F, \quad G,\quad H.$$
   The following table displays some useful information on these two-sided cells.
  \begin{center}
  \doublerulesep 0.4pt \tabcolsep8pt
\begin{tabular}{ccccc}
\hline
& &Number & Size of Jordan blocks& Maximal reductive subgroup\\
& & of left & of the  corresponding & of the centralizer of a unipotent \\
$X$ & {${a(X)}$} & cells in $X$ &unipotent class in $Sp_6(\mathbb C)$&element in the corresp. unipotent class\\
\hline$A$&9&48&\hfill (111111)\ \ \ \ \ \ \ \ \ \ \ \   &$Sp_6(\mathbb C)$\\
$B$&6&24&\hfill (21111)\ \ \ \ \ \ \ \ \ \ \ \  &$Sp_4(\mathbb C)\times\mathbb Z/2\mathbb Z$\\
$C$&4&18&\hfill (2211)\ \ \ \ \ \ \ \ \ \ \ \  &$SL_2(\mathbb C)\times O_2(\mathbb C)$\\
$D$&3&12&\hfill (222)\ \ \ \ \ \ \ \ \ \ \ \  &$O_3(\mathbb C)$\\
$E$&2&8&\hfill (411)\ \ \ \ \ \ \ \ \ \ \ \  &$SL_2(\mathbb C)\times\mathbb Z/2\mathbb Z$\\
$F$&2&6&\hfill (33)\ \ \ \ \ \ \ \ \ \ \ \  &$SL_2(\mathbb C)$\\
$G$&1&4&\hfill (42)\ \ \ \ \ \ \ \ \ \ \ \  &$\mathbb Z/2\mathbb Z\times\mathbb Z/2\mathbb Z$\\
$H$&0&1&\hfill (6)\ \ \ \ \ \ \ \ \ \ \ \  &$\mathbb Z/2\mathbb Z$\\
\hline
\end{tabular}
\end{center}\vspace{2mm}
The notations for two-sided cells in the table are the same as those in [D], which will be replaced by other notations in subsequent sections, otherwise confusion would happen since notations $C, F,G$ are already used for other objects.

In the rest of the paper, $W$ always stands for the {extended} affine Weyl group attached to $Sp_6(\mathbb C)$, $\tau,\ r_i$ are as in Subsection 2.1, and all representations are rational representations of algebraic groups.

\medskip

 \noindent{\bf 2.3.} We are concerned with the two-sided cell $c$ of $W$ containing $w_{012}=r_0r_1r_2r_0r_1r_2$. The corresponding unipotent class in $Sp_6(\mathbb C)$ has Jordan blocks (21111), and $u$ is a unipotent element in it.

Then $$F_c=Sp_4(\mathbb C)\times Z,$$
where $Z\simeq\mathbb Z/2\mathbb Z$ is the center of $Sp_6(\mathbb C)$. Since $Sp_4(\mathbb C)$ is simply connected, for any semisimple element $s\in F_c$ we have $Z_{F_c}(s)=Z_{F_c}(s)^\circ\times Z$ and $Z$ acts on $\mathcal B_u^s$ trivially. Therefore the set of isomorphism classes {of simple $\mathbf J_c$-modules} is one to one corresponding to the set of semisimple classes of $F_c$.

For a semisimple element $s$ in $F_c$, let $E_s$ be the corresponding simple $\mathbf J_c$-module. Then for a dominant weight $x\in X$, the element $\phi_c(S_x)$ acts on $E_s$ by scalar $\tr(ss_{q^{-1/2}}, V(x))$. This fact is useful in figuring out and verifying the map $\pi$ in Lusztig's conjecture stated in \S1.9.

 \section{The based ring of the two-sided cell containing $r_0r_1r_2r_0r_1r_2$}

\noindent{\bf 3.1.} In this section $c$ stands for the two-sided cell of $W$ containing $r_0r_1r_2r_0r_1r_2$.  And the value of $a$-function on $c$ is 6. According to [D, Figure I, Theorem 6.4], $c$ has 24 left cells.

 J. Du has given a representative for each left cell (see [D, Figure I, Theorem 6.4]. We recall his result for the two-sided cell $c$.

 (a) There are  24 left cells in the two-sided cell $c$ and a representative of each left cell in $ c$ are:
$$\begin{array}{ccccccc}
&\Gamma_{012},&012012; &\Gamma_{013},&0120123; &\Gamma_{2},&01201232;\\
&&&&&&\\
& \Gamma_{23},&012012323; &\Gamma_{02},&012012320;&\Gamma_{03},&0120123203; \\
&&&&&&\\
&\Gamma_{12},&012012321; &\Gamma_{13},&0120123213; &\Gamma_{01},&0120123201;\\
 &&&&&&\\
&\Gamma'_{013},&01201232013; &\Gamma'_{2},&012012320132;
 &\Gamma_{3},&0120123201323;\\
 &&&&&&\\
 &\Gamma'_{02}, &01201232032; &\Gamma'_{01}, &012012320321; &\Gamma'_{12}, &01201232132; \\
  &&&&&&\\
 &\widehat \Gamma'_{01}, &012012321320; 
 &\Gamma''_{12}, &0120123203212; &\Gamma'_{13}, &01201232032123; \\
 &&&&&&\\
 &\Gamma''_{2}, &012012320321232; 
 &\Gamma_{0}, &0120123203212320;&\Gamma''_{02}, &0120123213202; \\
 &&&&&&\\
  &\Gamma'_{03}, &01201232132023; 
 &\widehat \Gamma''_{2}, &012012321320232; &\Gamma_{1}, &0120123213202321,
 \end{array}$$
 
 here  we  simply write $i_1i_2\cdots i_k$ for a reduced expression  $r_{i_1}r_{i_2}\cdots r_{i_k}$ of an element in $W$ and the subscript  of the notation of each left cell indicates the $R$-set of the left cell, say the subscript 02 of $\Gamma'_{02}$ indicates $R(\Gamma'_{02})=\{r_0,\ r_2\}$.

 Let $\Gamma$ and $\Gamma'$ be two left cells of $W$. If $\Gamma'=\Gamma^*$ for some $*=\{r,t\}$, then we write $\Gamma\ \overset{\{r, t\}}{\text{------}}\ \Gamma'$. The following result is easy to verify.

{\bf Lemma 3.2.} Keep the notations above. Then we have four graphs:

 $$\begin{aligned} \Gamma_{012}&\ \overset{\{r_2, r_3\}}{\text{------}}\ \Gamma_{2}\ \overset{\{r_0, r_2\}}{\text{------}}\ \Gamma_{013}\\
 012012&\qquad 01201232\qquad 0120123\end{aligned}$$

 $$ \begin{aligned}\Gamma_{03}&\ \overset{\{r_0, r_2\}}{\text{------}}\ \Gamma_{23}\ \overset{\{r_1, r_2\}}{\text{------}}\ \Gamma_{13}\\
 0120123203&\qquad 012012323\qquad0120123213\end{aligned}$$

  $$ \begin{aligned}\Gamma'_2&\ \overset{\{r_0, r_2\}}{\text{------}}\ \Gamma'_{013}\ \overset{\{r_2, r_3\}}{\text{------}}\ \Gamma_{3}\\
  012012320132&\qquad 01201232013\qquad 0120123201323\end{aligned}$$

 $$\begin{aligned}
 \Gamma_{02}&\ \overset{\{r_2, r_3\}}{\text{------}}\ \Gamma'_{02}\ \overset{\{r_1, r_2\}}{\text{------}}\ \Gamma'_{01}\ \overset{\{r_0, r_2\}}{\text{------}}\ \Gamma''_{12}\ \overset{\{r_2, r_3\}}{\text{------}}\ \Gamma''_{2}\ \overset{\{r_0, r_2\}}{\text{------}}\ \Gamma_{0}\\
{\scriptstyle{012012320}} &\ {\scriptstyle{01201232032}}\ {\scriptstyle{012012320321}}\ {\scriptstyle{0120123203212}}\ {\scriptstyle{012012320321232}}\ {\scriptstyle{ 012012320321232}}\\
\ | & {\scriptstyle{\{r_1, r_2\}}}\qquad\qquad\qquad\qquad\qquad\qquad\qquad\qquad | {\scriptstyle{\{r_1, r_2\}}}\\
\Gamma_{01} &\qquad\qquad\qquad\qquad\qquad\qquad\qquad\qquad\qquad \Gamma'_{13}\\
{\scriptstyle{0120123201}} &\qquad \qquad\qquad\qquad\qquad\qquad\qquad\qquad {\scriptstyle{012010232032123}}\\
\ | & {\scriptstyle{\{r_0, r_2\}}}\\
\Gamma_{12} &\ \overset{\{r_2, r_3\}}{\text{------}}\ \Gamma'_{12}\ \overset{\{r_0, r_2\}}{\text{------}}\ \widehat \Gamma'_{01}\ \overset{\{r_1, r_2\}}{\text{------}}\ \Gamma''_{02}\ \overset{\{r_2, r_3\}}{\text{------}}\ \widehat \Gamma''_{2}\ \overset{\{r_1, r_2\}}{\text{------}}\ \Gamma_{1}\\
{\scriptstyle{012012321}} &\ {\scriptstyle{01201232132}}\ {\scriptstyle{012012321320}}\  {\scriptstyle{0120123213202}}\ {\scriptstyle{ 012012321320232}}\ {\scriptstyle{ 0120123213202321}}\\
&\qquad\qquad\qquad\qquad\qquad\qquad\qquad\qquad\qquad\qquad | {\scriptstyle{\{r_0, r_2\}}}\\
&\qquad\qquad\qquad\qquad\qquad\qquad\qquad\qquad\qquad\quad\ \ \Gamma'_{03}\\
&\qquad\qquad\qquad\qquad\qquad\qquad\qquad\quad\quad\quad{\scriptstyle{01201232132023}}
 \end{aligned}$$

 where  each vertex contains a left cell and  a representative element of the left cell.

\medskip

{\bf 3.3.}  Let
\begin{alignat*}{2}
Y_1&=\{\Gamma_{012}, \ \Gamma_{2},\  \Gamma_{013}\}; \\
Y_2&=\{\Gamma_{03}, \ \Gamma_{23}, \ \Gamma_{13}\}; \\
Y_3&=\{\Gamma'_{013},\ \Gamma'_{2},\ \Gamma_3\};\\
Y_4&=\{\Gamma_{01},\ \Gamma_{02},\ \Gamma'_{02}, \ \Gamma'_{01},\ \Gamma''_{12},\\
&\ \ \ \ \ \Gamma''_{2},\ \Gamma_{0},\ \Gamma'_{13},\ \Gamma_{12},\ \Gamma'_{12},\\
 &\ \ \ \ \ \hat \Gamma'_{01},\ \Gamma''_{02},\ \hat \Gamma''_{2},\ \Gamma_{1},\ \Gamma'_{03}\};\\
 Y&=Y_1\cup Y_2\cup Y_3\cup Y_4.\end{alignat*}

\medskip

 In the following we discuss the based ring of $c$. We know that  $F_c=Sp_4(\mathbb C)\times Z/2\mathbb Z$.  We only deal with rational  representations of $F_c$.

Let $W_1$ be the Weyl group of $Sp_4(\mathbb C)$ and let $\lambda_1,\ \lambda_2$ be fundamental weights of $Sp_4(\mathbb C)$. Then we have
\begin{align} W_1\cdot \lambda_1&=\{ \lambda_1,\ - \lambda_1+ \lambda_2,\  \lambda_1+ \lambda_2,\ - \lambda_1\}\\
W_1\cdot \lambda_2&=\{ \lambda_2, \ - 2\lambda_1+ \lambda_2,\  2\lambda_1-\lambda_2,\ - \lambda_2\}.\end{align}

 Let $F_c=Sp_4(\mathbb C)\times Z/2\mathbb Z$ act on $Y$ trivially. Then $K_{F_c}(Y\times Y)$ is isomorphic to the $Y\times Y$ matrix ring $M_Y(Rep\ F_c)$, where $Rep\ F_c$ is the representation ring of $F_c$.

  Let $\epsilon$ be the nontrivial one dimensional representation of $ \mathbb Z/2\mathbb Z$ and $V(\lambda)$ be an irreducible representation of $Sp_4(\mathbb C)$ with highest dominant weight $\lambda\in P^{+}:=\{a\lambda_1+b\lambda_2\ |\ a, b\in\mathbb N\}$. They can be regarded as irreducible representations of $F_c$ naturally. Up to isomorphism, the irreducible representations of $F_c$ are $V(\lambda),\ \epsilon\otimes V(\lambda), \lambda\in P^{+}$. We will denote $ \epsilon\otimes V(\lambda)$ by  $\epsilon  V(\lambda)$.

\bigskip

{\bf Theorem 3.4.}  Let $c$ be the two-sided cell of $W$ (the extended affine Weyl group attached to $Sp_6(\mathbb C)$) containing $w_{012}$, where $w_{012}$ is the longest element of the parabolic subgroup generated by $r_0, r_1, r_2$ and $\Gamma$ be any left cell in $c$. Then there exists a bijection
  $$\pi: c\to \text{the set of isomorphism classes of irreducible $F$-vector bundles on}\ Y\times Y,$$
  such that

  (a)  The bijection $\pi$ induces a based ring isomorphism
  $$\pi: J_c\to K_{F}(Y\times Y),\ \ t_x\mapsto \pi(x).$$

  (b) $\pi(x^{-1})_{(a,b)}=\pi(x)_{(b,a)}^*$ is the dual representation of $\pi(x)_{(b,a)}.$

 \medskip

  {\bf Remark.} This result actually has already been proved in Theorem 4 in [BO] conceptually since the group $F_c=Sp_4(\mathbb C)\times Z/2\mathbb Z$ has no nontrivial projective  representations. Here we will construct the bijection explicitly so that the isomorphism can be used to compute certain irreducible representations of affine Hecke algebras of type $\tilde B_3$.

  \medskip
\medskip

First we consider the based ring $J_{\Gamma\cap\Gamma^{-1}}$ for any left cell $\Gamma$ in $c$.

{\bf Theorem 3.4A.}
 There is a natural bijection
$$\pi_{\Gamma\cap \Gamma^{-1}}: \Gamma\cap \Gamma^{-1}\to \text{Irr}(Sp_4(\mathbb C)\times\mathbb Z/2\mathbb Z)$$
 and the bijection induces a based ring isomorphism
 $$\pi_{\Gamma\cap \Gamma^{-1}}:  J_{\Gamma\cap\Gamma^{-1}}\to \text{Rep}(Sp_4(\mathbb C)\times\mathbb Z/2\mathbb Z),\quad t_x\mapsto \pi(x),$$
 where Irr$(Sp_4(\mathbb C)\times\mathbb Z/2\mathbb Z)$ is the set of isomorphism classes of irreducible  representations of $Sp_4(\mathbb C)\times\mathbb Z/2\mathbb Z$, Rep$(Sp_4(\mathbb C)\times\mathbb Z/2\mathbb Z)$ is the Grothendieck group of the category of  representations of $Sp_4(\mathbb C)\times\mathbb Z/2\mathbb Z$, the multiplication in Rep$(Sp_4(\mathbb C)\times\mathbb Z/2\mathbb Z)$ is given by tensor product of representations.

{\bf Proof.} It is easy to see that if $\Gamma'=\Gamma^*$ for some $*=\{s,t\}$ with $s, t\in S$, then the map $\Gamma\cap\Gamma^{-1}\to\Gamma'\cap\Gamma'^{-1},\ x\mapsto{}^*x^*$ is a bijection. Hence by 1.4(g) we see that the map $J_{\Gamma\cap\Gamma^{-1}}\to J_{\Gamma'\cap\Gamma'^{-1}},\ t_x\mapsto t_{{}^*x^*}$ is a ring isomorphism.

By Lemma 3.2 ,1.4(f) and 1.4(g), it suffices to prove Theorem 3.4A  for $\Gamma=\Gamma_{012}, \Gamma_{03}, \Gamma'_{013}$ and $\Gamma_{02}$.
By Lemma 3.2, $\Gamma'_{013}$ can be obtained from $\Gamma_0$ by several $*$-operations. According to [B], for $\Gamma_0$, Theorem 3.4A  is true. By 1.4(g), Theorem 3.4A  is true for $\Gamma'_{013}$. Also, we offer here a different proof for $\Gamma'_{013}$ which exhibits the bijection $\pi_{\Gamma'_{013}\cap \Gamma_{013}^{-1}}: \Gamma'_{013}\cap \Gamma_{013}^{-1}\to \text{Irr}(Sp_4(\mathbb C)\times\mathbb Z/2\mathbb Z)$ explicitly.

And we  also need to consider the $t_{d_\Theta}Jt_{d_\Theta}-t_{d_\Gamma}Jt_{d_\Gamma}$-bimodule $J_{\Gamma\cap\Theta^{-1}}$ when left cells $\Gamma$ and $\Theta$ are not in the same $Y_i$ for any $1\le i\le 4$.

{\bf Theorem 3.4B:} The bijection
$$\pi_{ \Gamma\cap \Theta^{-1}}:  \Gamma\cap \Theta^{-1}\to\text{the set of isomorphism classes of irreducible $F_c$-vector bundles on}\  \Theta\times \Gamma,$$
  induces an isomorphism of $\mathbb Z$-modules
 $$\pi_{\Gamma\cap \Theta^{-1}}: J_{\Gamma\cap \Theta^{-1}}\to K_{F_c}(\Theta\times \Gamma)\cong\text{Rep}\ (Sp_4(\mathbb C)\times\mathbb Z/2\mathbb Z),$$
   such that take any $x\in\Gamma\cap \Theta^{-1}$, $\Phi\in Y_1\cup Y_2\cup Y_3\cup Y_4$, we have

(a) $\pi_{\Phi\cap\Theta^{-1}}(t_xt_y)=\pi_{\Gamma\cap\Theta^{-1}}(t_x)\pi_{\Phi\cap\Gamma^{-1}}(t_y), \text{for\ any}\ y\in \Phi\cap\Gamma^{-1}$,

(b) $\pi_{\Gamma\cap\Phi^{-1}}(t_zt_x)=\pi_{\Theta\cap\Phi^{-1}}(t_z)\pi_{\Gamma\cap\Theta^{-1}}(t_x), \text{for\ any}\ z\in \Theta\cap\Phi^{-1}$.

{\bf Proof.} It is easy to see that if $\Gamma'=\Gamma^*, \Theta'=\Theta^\sharp$ for some $*=\{r,t\}, \sharp=\{r',t'\}$ with $r, t, r', t'\in S$, then the map $\Gamma\cap\Theta^{-1}\to\Gamma'\cap\Theta^{-1},\ x\mapsto{}^\sharp x^*$ is a bijection. Hence by 1.4(g) we see that the map $J_{\Gamma\cap\Theta^{-1}}\to J_{\Gamma'\cap\Theta'^{-1}},\ t_x\mapsto t_{{}^\sharp x^*}$ is a $\mathbb Z$-module isomorphism.

By Lemma 3.2 ,1.4(e), 1.4(f) and 1.4(g), it suffices to prove Theorem 3.4B  for $\Gamma\cap\Theta^{-1}=\Gamma_{012}\cap\Gamma_{02}^{-1}, \Gamma_{012}\cap\Gamma_{03}^{-1}, \Gamma_{012}\cap(\Gamma'_{013})^{-1}, \Gamma_{02}\cap\Gamma_{03}^{-1}, \Gamma_{02}\cap(\Gamma'_{013})^{-1}, \Gamma_{03}\cap(\Gamma'_{013})^{-1}$.

  The rest of this section is devoted to proving the theorems above.

     Note that for any $x$ in $ \Gamma\cap \Gamma^{-1}$ with any left cell $\Gamma$ in $c$, we have $x=x^{-1}$. Using 1.4(e) we get

    (c) $t_xt_y=t_yt_x$ for any $x,y\in \Gamma\cap \Gamma^{-1}$.

    (d) Denote $t_xt_y=\sum_{z\in W} \gamma_{x,y,z}t_z$ for any $x,y\in \Gamma\cap \Gamma^{-1}$. Then $t_{\tau x}t_y=\sum_{z\in W} \gamma_{x,y,z}t_{\tau z}$.

\medskip

  {\bf 3.5. }   In this subsection we consider $\Gamma_{012}\cap\Gamma_{012}^{-1}$ and its based ring $J_{\Gamma_{012}\cap\Gamma_{012}^{-1}}$.

   Let $*=\{r_0, r_2\}$, $\sharp=\{r_1, r_2\}$ and $\star=\{r_2, r_3\}$. According to  [D, Theorem 6.4],  the set
$$ \Gamma_{012}\cap \Gamma_{012}^{-1}=\{ w_{012}(\tau r_3r_2r_0r_1r_2)^i(r_3r_2r_0r_1)^{2j}, \tau w_{012}(\tau r_3r_2r_0r_1r_2)^i(r_3r_2r_0r_1)^{2j},\,|\, i, j\ge 0\}.$$

  {\bf Proposition 3.5A.}
  Keep the notations from subsection 3.1, 3.2 and 3.3. That is,
 $c$ stands for the two-sided cell of $W$  containing $w_{012}$, and $\Gamma_{012}$ be the left cell in $c$ containing $w_{012}$. Then the map
\begin{alignat*}{2} \pi_{\Gamma_{012}\cap  \Gamma_{012}^{-1}}:    \Gamma_{012}\cap  \Gamma_{012}^{-1}&\longrightarrow \text{Irr}(Sp_4(\mathbb C)\times\mathbb Z/2\mathbb Z),\\
w_{012}(\tau r_3r_2r_0r_1r_2)^i(r_3r_2r_0r_1)^{2j}&\longmapsto V(i\lambda_1+j\lambda_2),\\
\tau w_{012}(\tau r_3r_2r_0r_1r_2)^i(r_3r_2r_0r_1)^{2j}&\longmapsto \epsilon V(i\lambda_1+j\lambda_2),
\end{alignat*}
induces a based ring isomorphism
    $$\pi_{\Gamma_{012}\cap \Gamma_{012}^{-1}}:  J_{\Gamma_{012}\cap \Gamma_{012}^{-1}}\longrightarrow \text{Rep}(Sp_4(\mathbb C)\times\mathbb Z/2\mathbb Z),\quad t_x\longmapsto  \pi_{\Gamma_{012}\cap  \Gamma_{012}^{-1}}(x).$$

     \medskip

  Let $w_0$ be  the longest element of the Weyl group $W_0$ and $X^+=\{x\in X | l(w_0x)=l(w_0)+l(x)\}$ be the set of dominant weights in $X$. Denote by $x_1,x_2,\ x_3$ be the fundamental weights in $X^+$ corresponding to simple reflections $r_1,\ r_2,\ r_3$ respectively.  Then we have
  \begin{equation}
   S_{x_1}=\theta_{x_1}+\theta_{x^{-1}_1x_2}+\theta_{x^{-1}_2x_3}+
 \theta_{x^{-1}_1}+\theta_{x_1x^{-1}_2}+\theta_{x_2x^{-1}_3},
  \end{equation}
  \begin{equation}
  \begin{aligned}
  S_{x_2}
  =&\theta_{x_2}+\theta_{x_1x^{-1}_2x_3}+\theta_{x_1x_2x^{-1}_3}+
  \theta_{x^2_1x^{-1}_2}+\theta_{x^{-1}_1x_3}+\theta_{x^{-1}_1x^2_2x^{-1}_3}\\
  &+\theta_{x^{-1}_2}+\theta_{x_1^{-1}x_2x^{-1}_3}+\theta_{x^{-1}_1x^{-1}_2x_3}
  +\theta_{x^{-2}_1x_2}+\theta_{x_1x^{-1}_3}+\theta_{x_1x^{-2}_2x_3}
+2\tilde T_e,
\end{aligned}
  \end{equation}

  where $e$ is the identity element in $W$. 

    Set
    \begin{equation}x_{i,j}=w_{012}(\tau r_3r_2r_0r_1r_2)^i(r_3r_2r_0r_1)^{2j},\quad \forall i, j\geq 0.\end{equation}

    Before continuing, we make a convention: {\sl we shall use the symbol $\Box$ for any element in the two-sided ideal $H^{<w_{012}}$ of $H$ spanned by all $C_u$ with $a(u)>6$. Then $\Box+\Box=\Box$ and $h\Box=\Box$ for any $h\in H$.}

    First we need to verify the following Lemma 3.5B and Lemma 3.5C.

{\bf Lemma 3.5B.}
Set $t_{x_{-1, j}}=t_{x_{i, -1}}=0$. Then for any $i, j\ge 0$, we have
$$t_{x_{1,0}}t_{x_{i,j}}=
    t_{x_{i+1,j}}+t_{x_{i+1,j-1}}+t_{x_{i-1,j+1}}+t_{x_{i-1,j}}.$$

\begin{proof}

We first compute $S_{x_1}C_{x_{i,j}}$.

By (3),  we have \begin{equation}S_{x_1}C_{x_{i,j}}=(\theta_{x_1}+\theta_{x^{-1}_1x_2}+\theta_{x^{-1}_2x_3}+\theta_{x^{-1}_1}
+\theta_{x_1x^{-1}_2}+\theta_{x_2x^{-1}_3})C_{x_{i,j}}.\end{equation}

Note that by 1.1(b) for any $r\in S$ and $u\in W$, we have

 $$\tilde T_rC_u=\begin{cases}\displaystyle q^{\frac12}C_u,\quad &\text{if\ }ru\leq u,\\
\displaystyle C_rC_u-q^{-\frac12}C_u,\quad &\text{if\ }ru>u\end{cases}.$$
$$\tilde T_r^{-1}C_u=\begin{cases}\displaystyle q^{-\frac12}C_u,\quad &\text{if\ }ru\leq u,\\
\displaystyle C_rC_u-q^{\frac12}C_u,\quad &\text{if\ }ru>u\end{cases}.$$

According to 1.7(c), for any $y\in\Gamma_{012}\cap\Gamma_{012}^{-1},$ we have
\begin{equation}S_{x_1}C_y=\sum_{z\in \Gamma_{012}\cap\Gamma_{012}^{-1}}\zeta_{x_1,y,z}C_z+\Box.\end{equation}

We shall use $\bigtriangleup$ for any linear combination of $C_z,\ z\not\in\Gamma_{012}\cap\Gamma_{012}^{-1}$. Then $\bigtriangleup+\bigtriangleup=\bigtriangleup$.

Set $C_{x_{-1,0}}=C_{x_{0,-1}}=0$ and denote $\tilde T_{r_i}$ by $\tilde T_i$ by simplicity. By page 644 in [L1] and direct computation we get

\begin{itemize}
\item $\theta_{x_1}C_{x_{i,j}}=\tilde T_\tau\tilde T_1\tilde T_2\tilde T_3\tilde T_2\tilde T_1C_{x_{i,j}}=\bigtriangleup$.
\item $\theta_{x^{-1}_1x_2}C_{x_{i,j}}=\tilde T^{-1}_1\tilde T_{x_1}\tilde T^{-1}_1C_{x_{i,j}}=\tilde T_\tau\tilde T^{-1}_{0}\tilde T_1\tilde T_2\tilde T_3\tilde T_2C_{x_{i,j}}=\bigtriangleup.$
\vspace*{-1em}
\begin{align*}
\bullet \,\, \theta_{x^{-1}_2x_3}C_{x_{i,j}}=&\tilde T^{-1}_2\tilde T^{-1}_1\tilde T_{x_1}\tilde T^{-1}_1\tilde T^{-1}_2C_{x_{i,j}}
=\tilde T_\tau\tilde T^{-1}_{2}\tilde T^{-1}_{0}\tilde T_1\tilde T_2\tilde T_3C_{x_{i,j}}\qquad\qquad\qquad\\
=&C_{x_{i+1,j}}+C_{x_{i-1,j+1}}+C_{x_{i-1,j}}+C_{x_{i+1,j-1}}- q^{-\frac12}C_{\tau x_{i,j}}+\bigtriangleup.\end{align*}
\vspace*{-2em}
\item $\theta_{x_1^{-1}}C_{x_{i,j}}=\tilde T_\tau \tilde T^{-1}_0\tilde T^{-1}_2\tilde T^{-1}_3\tilde T^{-1}_2\tilde T^{-1}_0C_{x_{i,j}}=\bigtriangleup.$
\item $\theta_{x_1x^{-1}_2}C_{x_{i,j}}=\tilde T_1\tilde T^{-1}_{x_1}\tilde T_1C_{x_{i,j}}=\tilde T_\tau\tilde T^{-1}_{2}\tilde T^{-1}_3\tilde T^{-1}_2\tilde T^{-1}_0\tilde T_1C_{x_{i,j}}=\bigtriangleup.$
\item $\theta_{x_2x^{-1}_3}C_{x_{i,j}}=\tilde T_2\tilde T_1\tilde T^{-1}_{x_1}\tilde T_1\tilde T_2C_{x_{i,j}}=\tilde T_\tau\tilde T^{-1}_3\tilde T^{-1}_2\tilde T^{-1}_0\tilde T_1\tilde T_2C_{x_{i,j}}= -q^{\frac12}C_{\tau x_{i,j}}+\bigtriangleup.$

\end{itemize}

Combining the above six identities and formulas (6)-(7) we get
\begin{align} S_{x_1}C_{x_{i,j}}=C_{x_{i+1,j}}+C_{x_{i-1,j+1}}+C_{x_{i-1,j}}+C_{x_{i+1,j-1}}- (q^{\frac12}+q^{-\frac12})C_{\tau x_{i,j}}+\Box.
\end{align}
When $i=j=0$, $x_{i,j}=w_{012}$. Let $\xi=q^{\frac12}+q^{-\frac12}.$ Since $C_{x_{-1,0}}=C_{x_{0,-1}}=0,$ we get
\begin{align} S_{x_1}C_{w_{012}}=C_{x_{1,0}}-\xi C_{\tau w_{012}}+\Box.
\end{align}

Let $C_{w_{012}}C_{w_{012}}=\eta C_{w_{012}}$. Then $\eta=q^{-3}\sum_{u\in W_{012}}q^{l(u)},$ here $W_{012}$ is the parabolic subgroup of $W$ generated by $r_0,\ r_1,\ r_2$. Clearly we have
\begin{align}C_{w_{012}}C_{x_{i,j}}=\eta C_{x_{i,j}},\quad C_{\tau w_{012}}C_{x_{i,j}}=\eta C_{\tau x_{i,j}}.\end{align}

 Using formulas (8)-(10), we get
\begin{equation}\begin{array}{ll}C_{x_{1,0}}C_{x_{i,j}}&=(S_{x_1}C_{w_{012}}+\xi C_{\tau w_{012}})C_{x_{i,j}}+\Box\\
&=\eta S_{x_1}C_{x_{i,j}}+\xi\eta C_{\tau x_{i,j}}+\Box\\
&=\eta(C_{x_{i+1,j}}+C_{x_{i-1,j+1}}+C_{x_{i-1,j}}+C_{x_{i+1,j-1}}- \xi C_{\tau x_{i,j}}+\Box)+\xi\eta C_{\tau x_{i,j}}+\Box\\
&=\eta(C_{x_{i+1,j}}+C_{x_{i-1,j+1}}+C_{x_{i-1,j}}+C_{x_{i+1,j-1}})+\Box.\end{array}\end{equation}
Hence we have
 $$t_{x_{1,0}}t_{x_{i,j}}=t_{x_{i+1,j}}+t_{x_{i-1,j+1}}+t_{x_{i-1,j}}+t_{x_{i+1,j-1}}.$$
  The proof is completed.
\end{proof}

{\bf Lemma 3.5C.}
Set $t_{x_{-1, j}}=t_{x_{i, -1}}=0$. Then for any $i, j\ge 0$, we have
$$t_{x_{0,1}}t_{x_{i, j}}=
    t_{x_{i, j+1}}+t_{x_{i+2, j-1}}+(1-\delta_{0,i})t_{x_{i, j}}+t_{x_{i-2,j+1}}+t_{x_{i, j-1}}.$$

\begin{proof}
First of all, the case when $i=1, j=0$ follows from Lemma 3.5B.

Let $\xi=q^{\frac12}+q^{-\frac12}.$
Next, we compute $S_{x_2}C_{w_{012}}$.

According to 1.7(c), for any $y\in\Gamma_{012}\cap\Gamma_{012}^{-1},$ we have
\begin{equation}S_{x_2}C_y=\sum_{z\in \Gamma_{012}\cap\Gamma_{012}^{-1}}\zeta_{x_2,y,z}C_z+\Box.\end{equation}

Similar to the proof of Lemma 3.5B, we use $\bigtriangleup$ for any linear combination of $C_z,\ z\not\in\Gamma_{012}\cap\Gamma_{012}^{-1}$, set $C_{x_{-1,0}}=C_{x_{0,-1}}=0$ and denote $\tilde T_{r_i}$ by $\tilde T_i$ by simplicity. By page 644 in [L1], formula (4) and direct computation we get

\begin{itemize}
\item $\theta_{x_2}C_{x_{i,j}}=\tilde T_0\tilde T_2\tilde T_3\tilde T_2\tilde T_1\tilde T_2\tilde T_3\tilde T_2C_{x_{i,j}}=\bigtriangleup$.
\item $\theta_{x_1x^{-1}_2x_3}C_{x_{i,j}}=\tilde T_0\tilde T_2\tilde T_3\tilde T^{-1}_0\tilde T_2\tilde T_1\tilde T_2\tilde T_3C_{x_{i,j}}=\bigtriangleup.$
\item $\theta_{x_1x_2x^{-1}_3}C_{x_{i,j}}=\tilde T_0\tilde T^{-1}_3\tilde T_2\tilde T_3\tilde T^{-1}_{0}\tilde T_{2}\tilde T_1\tilde T_2C_{x_{i,j}}=\bigtriangleup.$
\item $\theta_{x^2_1x^{-1}_2}C_{x_{i,j}}=\tilde T^{-1}_2\tilde T_0\tilde T^{-1}_3\tilde T_2\tilde T_3\tilde T^{-1}_0\tilde T_2\tilde T_1C_{x_{i,j}}=-C_{x_{i,j}}+\bigtriangleup.$
\vspace*{-1em}
\begin{align*}
\bullet \, \, \theta_{x^{-1}_1x_3}C_{x_{i,j}}=&\tilde T^{-1}_1\tilde T^{-1}_{2}\tilde T_0\tilde T_2\tilde T_1\tilde T_3\tilde T_2\tilde T_3C_{x_{i,j}}\\
=&C_{x_{i,j+1}}+C_{x_{i+2,j-1}}+(1-\delta_{0,i})C_{x_{i,j}}+C_{x_{i,j-1}}+C_{x_{i-2,j+1}}\qquad\qquad\\
&-q^{-\frac{1}{2}}(C_{\tau x_{i+1,j}}+C_{\tau x_{i+1,j-1}}+C_{\tau x_{i-1,j}}+C_{\tau x_{i-1,j+1}})+\bigtriangleup.\end{align*}
\vspace*{-2.5em}
\item $\theta_{x^{-1}_1x^2_2x^{-1}_3}C_{x_{i,j}}=\tilde T^{-1}_3\tilde T_0\tilde T_2\tilde T_1\tilde T^{-1}_2\tilde T_3\tilde T^{-1}_0\tilde T_2C_{x_{i,j}}=\bigtriangleup.$
\item $\theta_{x^{-1}_2}C_{x_{i,j}}=\tilde T^{-1}_2\tilde T^{-1}_3\tilde T^{-1}_2\tilde T^{-1}_1\tilde T^{-1}_2\tilde T^{-1}_3\tilde T^{-1}_2\tilde T^{-1}_0C_{x_{i,j}}=-q^{-1}C_{x_{i,j}}+\bigtriangleup.$
\item $\theta_{x_1^{-1}x_2x^{-1}_3}C_{x_{i,j}}=\tilde T^{-1}_3\tilde T^{-1}_2\tilde T^{-1}_1\tilde T^{-1}_2\tilde T_0\tilde T^{-1}_3\tilde T^{-1}_2\tilde T^{-1}_0C_{x_{i,j}}=\bigtriangleup.$
\item $\theta_{x^{-1}_1x^{-1}_2x_3}C_{x_{i,j}}=\tilde T^{-1}_2\tilde T^{-1}_1\tilde T^{-1}_2\tilde T_0\tilde T^{-1}_3\tilde T^{-1}_2\tilde T_3\tilde T^{-1}_0C_{x_{i,j}}=(q^{-1}-1)C_{x_{i,j}}+\bigtriangleup.$
\item $\theta_{x^{-2}_1x_2}C_{x_{i,j}}=\tilde T^{-1}_1\tilde T^{-1}_2\tilde T_0\tilde T^{-1}_3\tilde T^{-1}_2\tilde T_3\tilde T^{-1}_0\tilde T_2C_{x_{i,j}}=\bigtriangleup.$
\item $\theta_{x_1x^{-1}_3}C_{x_{i,j}}=\tilde T^{-1}_3\tilde T^{-1}_2\tilde T^{-1}_3\tilde T^{-1}_1\tilde T^{-1}_2\tilde T^{-1}_0\tilde T_2\tilde T_1C_{x_{i,j}}=\bigtriangleup.$\vspace*{-1em}
 \begin{align*}
\bullet \, \, \theta_{x_1x^{-2}_2x_3}C_{x_{i,j}}&=\tilde T^{-1}_2\tilde T_0\tilde T^{-1}_3\tilde T_2\tilde T^{-1}_1\tilde T^{-1}_2\tilde T^{-1}_0\tilde T_3C_{w_{012}}\\&=-q^{\frac{1}{2}}(C_{\tau x_{i+1,j}}+C_{\tau x_{i+1,j-1}}+C_{\tau x_{i-1,j}}+C_{\tau x_{i-1,j+1}})+C_{x_{i,j}}+\bigtriangleup.\end{align*}
\end{itemize}

Combining the above 12 identities, formula (4) and (12), we get
\begin{equation}\begin{aligned}S_{x_2}C_{x_{i,j}}=&C_{x_{i,j+1}}+C_{x_{i+2,j-1}}+C_{x_{i,j-1}}+C_{x_{i-2,j+1}}-\xi (C_{\tau x_{i+1,j}}+C_{\tau x_{i+1,j-1}}+C_{\tau x_{i-1,j}}+C_{\tau x_{i-1,j+1}})\\
&+(2-\delta_{0,i})C_{x_{i,j}}+\Box.\end{aligned}\end{equation}
\begin{align}S_{x_2}C_{w_{012}}=C_{x_{0,1}}-\xi C_{\tau x_{1,0}}+C_{w_{012}}+\Box.\end{align}

Thus by formula (10) and (14), we get \begin{align}C_{x_{0,1}}C_{x_{i,j}}=(S_{x_2}C_{w_{012}}+\xi C_{\tau x_{1,0}}-C_{w_{012}})C_{x_{i,j}}+\Box=\eta S_{x_2}C_{x_{i,j}}+\xi C_{\tau x_{1,0}}C_{x_{i,j}}-\eta C_{x_{i,j}}+\Box.\end{align}

Combining it with (11) and (13), we get
\begin{equation}\begin{aligned}
C_{x_{0,1}}C_{x_{i,j}}
=&\eta (C_{x_{i,j+1}}+C_{x_{i+2,j-1}}+C_{x_{i,j-1}}+C_{x_{i-2,j+1}}+(1-\delta_{0,i})C_{x_{i,j}})+\Box\\
\in& \xi^6(C_{x_{i,j+1}}+C_{x_{i+2,j-1}}+C_{x_{i,j-1}}+C_{x_{i-2,j+1}}+(1-\delta_{0,i})C_{x_{i,j}})+H^{<w_{012}},
\end{aligned}\end{equation}

Therefore, $t_{x_{0,1}}t_{x_{i, j}}=
    t_{x_{i, j+1}}+t_{x_{i+2, j-1}}+(1-\delta_{0,i})t_{x_{i, j}}+t_{x_{i-2,j+1}}+t_{x_{i, j-1}}.$
Then the proof is completed.
\end{proof}

Set $t_{x_{m, n}}=0$ if $m<0$ or $n<0$.
Then for any $i,j,k,l\ge 0$, we have
$$t_{x_{k+1,l}}t_{x_{i,j}}=(t_{x_{1,0}}t_{x_{k,l}}+t_{x_{1,0}}t_{x_{k-2,l}}-t_{x_{0,1}}t_{x_{k-1,l}}-\delta_{0,k-1}t_{x_{k-1,l}}-t_{x_{k-1,l-1}}-t_{x_{k-3,l}})t_{x_{i,j}},$$
$$t_{x_{k,l+1}}t_{x_{i,j}}=(t_{x_{0,1}}t_{x_{k,l}}+t_{x_{0,1}}t_{x_{k,l-1}}-t_{x_{1,0}}t_{x_{k-1,l}}-t_{x_{1,0}}t_{x_{k+1,l-1}}-t_{x_{k,l-2}})t_{x_{i,j}},$$
So by induction, we can get the product of any elements in $J_{\Gamma_{012}\cap\Gamma_{012}^{-1}}$.

On the other hand, in $\text{Rep}(Sp_4(\mathbb C)\times \mathbb Z/2\mathbb Z)$, we have
$$V(\lambda_1)V(i\lambda_1+j\lambda_2)=V(\lambda_{i+1}+\lambda_j)+V(\lambda_{i+1}+\lambda_{j-1})+V(\lambda_{i-1}+\lambda_{j+1})+V(\lambda_{i-1}+\lambda_{j}),$$
$$V(\lambda_2)V(i\lambda_1+j\lambda_2)=V(\lambda_{i}+\lambda_{j+1})+V(\lambda_{i+2}+\lambda_{j-1})+(1-\delta_{0,i})V(\lambda_{i}+\lambda_{j})+V(\lambda_{i-2}+\lambda_{j+1})+V(\lambda_{i}+\lambda_{j-1}).$$

{\bf Corollary 3.5D.} For any left cell $\Gamma,\ \Theta$ in $Y_1$,  we have a bijection
  $$\pi_{\Gamma \cap\Theta^{-1}}: \Gamma \cap\Theta^{-1}\to\text{Irr}\ (Sp_4(\mathbb C)\times\mathbb Z/2\mathbb Z),$$

  inducing an isomorphism of $\mathbb Z$-modules
  $$\pi_{\Gamma\cap\Theta^{-1}}: J_{\Gamma \cap\Theta^{-1}}\to\text{Rep}\ (Sp_4(\mathbb C)\times\mathbb Z/2\mathbb Z),$$
such that take any $x\in\Gamma \cap\Theta^{-1}$, $\Phi\in Y_1$, we have

(a) $\pi_{\Phi\cap\Theta^{-1}}(t_xt_y)=\pi_{\Gamma\cap\Theta^{-1}}(t_x)\pi_{\Phi\cap\Gamma^{-1}}(t_y), \text{for\ any}\ y\in \Phi\cap\Gamma^{-1}$,

(b) $\pi_{\Gamma\cap\Phi^{-1}}(t_zt_x)=\pi_{\Theta\cap\Phi^{-1}}(t_z)\pi_{\Gamma\cap\Theta^{-1}}(t_x), \text{for\ any}\ z\in \Theta\cap\Phi^{-1}$.

\begin{proof}
It follows Lemma 3.2 and subsection 3.3.
\end{proof}

   \medskip

 \noindent{\bf 3.6.} In this subsection we establish a bijection

  $$\pi_{\Gamma_{02}\cap \Gamma_{02}^{-1}}: \Gamma_{02}\cap\Gamma_{02}^{-1}\to\text{Irr}(Sp_4(\mathbb C)\times\mathbb Z/2\mathbb Z)$$
  which induces an isomorphism of based rings
  $$\pi_{\Gamma_{02}\cap \Gamma_{02}^{-1}}: J_{\Gamma_{02}\cap \Gamma^{-1}_{02}}\to \text{Rep}(Sp_4(\mathbb C)\times\mathbb Z/2\mathbb Z).$$

Recall the notations in subsection 3.5 and (5), then
 according to [D, Theorem 6.4] we have
$$ \begin{aligned}
&\Gamma_{02}\cap \Gamma_{02}^{-1}\\
=&\{ r_0r_2r_3w_{012}(\tau r_3r_2r_0r_1r_2)^i(r_3r_2r_0r_1)^{2j}r_3r_2r_0, \tau r_1r_2r_3w_{012}(\tau r_3r_2r_0r_1r_2)^i(r_3r_2r_0r_1)^{2j}r_3r_2r_0 \,|\, i, j\ge 0\}\\
=&\{r_0({}^\star x_{i,j}^\star)r_0, r_0({}^\star \tau x_{i,j}^\star)r_0\,|\, x_{i,j}\in\Gamma_{012}\cap\Gamma_{012}^{-1}, i,j\ge 0\}\\
=&\{r_0({}^\star x{}^\star)r_0\,|\, x\in\Gamma_{012}\cap\Gamma_{012}^{-1}\}.\end{aligned}$$

{\bf Proposition 3.6A.} The bijection
\begin{alignat*}{2} \pi_{\Gamma_{02}\cap \Gamma_{02}^{-1}}:    \Gamma_{02}\cap  \Gamma_{02}^{-1}&\longrightarrow \text{Irr}(Sp_4(\mathbb C)\times\mathbb Z/2\mathbb Z),\\
r_0xr_0&\longmapsto \pi_{\Gamma_{012}\cap  \Gamma_{012}^{-1}}(x)
\end{alignat*}
   induces a based ring isomorphism
    $$ \pi_{\Gamma_{02}\cap \Gamma_{02}^{-1}}:  J_{\Gamma_{02}\cap \Gamma_{02}^{-1}}\longrightarrow \text{Rep}(Sp_4(\mathbb C)\times\mathbb Z/2\mathbb Z),\quad t_y\longmapsto \pi_{\Gamma_{02}\cap  \Gamma_{02}^{-1}}(y).$$

    \medskip

    Let
    $$\begin{aligned}\bar x_{i,j}={}^\star x_{i,j}^\star=r_2r_3x_{i,j}r_3r_2, \forall i, j\geq 0.\end{aligned}$$ Then $\Gamma_{02}\cap \Gamma_{02}^{-1}=\{r_0\bar x_{i,j}r_0 | \forall i, j\geq 0\}$.
    Thanks to identity (c) in subsection 3.4 and the discussion in subsection 3.5, to see Proposition 3.6A we only need to verify the following identity: for any $i,j,k,l\ge 0$, $z\in W$,
\begin{equation}\gamma_{r_0\bar x_{i,j}r_0, r_0\bar x_{k,l}r_0, r_0(^{\star}z^{\star})r_0}=\gamma_{\bar x_{i,j}, \bar x_{k,l}, ^{\star}z^{\star}}=\gamma_{x_{i,j}, x_{k,l}, z}.\end{equation}

\begin{proof}

Throughout the proof, we denote $x_{i,j}=w_{012}x'$ for $x'\in W$ with $l(x_{i,j})=l(x')+6$.

{\bf Step 1:} Compute $C_{r_0}C_{r_2r_3x_{i,j}v}$, for $v\in\{e,r_3,r_3r_2,r_3r_2r_0\}$.

By 1.1(b),  we have $C_{r_0}C_{r_2r_3x_{i,j}v}=C_{r_0r_2r_3x_{i,j}v}+\sum\limits_{\substack{y\prec r_2r_3x_{i,j}v\\ r_0y<y}}\mu(y, r_2r_3x_{i,j}v)C_y.$
 Note that ${L}(r_2r_3x_{i,j}v)=\{r_2\}$.

Assume $y\prec r_2r_3x_{i,j}v$ and $r_0y<y$.
\begin{itemize}
\item If $r_2y>y$, then by 1.1(a) we get $r_2y=r_2r_3x_{i,j}v$. So $y=r_3x_{i,j}v$.
\item If $r_0y<y$ and $r_2y<y$, then $y=r_0r_2r_0y'$ for some $y'\in W$ with $l(y)=l(y')+3$.
Note if $r_1\in{L}(y)$, then $a(y)\ge a(w_{012})=6$ and $C_y\in H^{<w_{012}}$. Next assume $r_1\not\in{L}(y)$. So by 1.3(d) we get $\mu(y, r_2r_3x_{i,j}v)=\tilde\mu({}^\sharp y, {}^\sharp(r_2r_3x_{i,j}v))=\tilde\mu({}^\sharp y, r_3x_{i,j}v)$ with ${L}(r_3x_{i,j}v)=\{r_0,r_1,r_3\}$.
\begin{itemize}
\item If ${}^\sharp y=r_0r_2y'$, then $y=r_2r_0r_2r_1r_2\tilde y $ for some $\tilde y\in W$ with $l(y)=l(\tilde y)+5$ and $\mu(y, r_2r_3x_{i,j}v)=\mu(r_0r_2r_1r_2\tilde y, r_3x_{i,j}v)$. If $r_3\in{L}(r_0r_2r_1r_2\tilde y)$, then $a(y)\ge a(r_0r_2r_1r_2\tilde y)\ge a(r_2r_1r_2\tilde y)=a(w_{123})=9$ and $C_y\in H^{<w_{012}}$. Otherwise by 1.1(a) we get $r_3r_0r_2r_1r_2\tilde y=r_3x_{i,j}v$. In this case, $y=r_2x_{i,j}v$, contradicting $r_2y<y$. So ${}^\sharp y=r_0r_2y'$ would not occur.
\item If ${}^\sharp y=r_1y$, then $\mu(y, r_2r_3x_{i,j}v)=\tilde\mu(r_1y, r_3x_{i,j}v)$.
If $l(y)=l(r_2r_3x_{i,j}v)-1$, then $\mu(y, r_2r_3x_{i,j}v)=\mu(r_3x_{i,j}v, r_1y)$ with $l(r_3x_{i,j}v)=l(r_1y)-1$. 
So $y\prec r_2r_3x_{i,j}v$ forces $r_3r_2r_0r_1r_2x'v\prec r_2r_0y'$. Since $r_2\not\in{L}(r_3r_2r_0r_1r_2x'v)$, by 1.1(a) we get $r_2r_3r_2r_0r_1r_2x'v=r_2r_0y'$. Thus $y=r_0r_2r_3r_2r_0r_1r_2x'v$, contradicting $r_2y<y$.
If $l(y)\le l(r_2r_3x_{i,j}v)-3$, then $\mu(y, r_2r_3x_{i,j}v)=\mu(r_1y, r_3x_{i,j}v)$. If $r_3\in{L}(r_1y)$, then $a(y)\ge a(w_{023})=9$ and $C_y\in H^{<w_{012}}$. Otherwise by 1.1(a) we get $r_3r_1y=r_3x_{i,j}v$ and $y=r_0r_2r_0r_1r_2x'v$.
This contradicts  ${}^\sharp y=r_1y$.

In conclusion, ${}^\sharp y=r_1y$ only deduces $C_y\in H^{<w_{012}}$.
\end{itemize}\end{itemize}

Therefore,
\begin{align}
C_{r_0}C_{r_2r_3x_{i,j}v}=C_{r_0r_2r_3x_{i,j}v}+C_{r_3x_{i,j}v}+\Box,\end{align}
By 1.4(e) we get
\begin{align}C_{v^{-1}x_{i,j}r_3r_2}C_{r_0}=C_{v^{-1}x_{i,j}r_3r_2r_0}+C_{v^{-1}x_{i,j}r_3}+\Box.\end{align}

Note for $u\in\{e, r_3, r_3r_2\}$,
$s_2s_3x_{i,j}u=\begin{cases}{^{\star}x_{i,j}}&{\text{if}\ u=e}\\
{^{\star}x_{i,j}^{\star*}}&{\text{if}\ u=r_3}\\
{^{\star}x_{i,j}^{\star}}&{\text{if}\ u=r_3r_2}\end{cases}$,\quad
$u^{-1}x_{i,j}s_3s_2=\begin{cases}{x_{i,j}^{\star}}&{\text{if}\ u=e}\\
{^{*\star}x_{i,j}^{\star}}&{\text{if}\ u=r_3}\\
{^{\star}x_{i,j}^{\star}}&{\text{if}\ u=r_3r_2}\end{cases}.$

Hence by 1.4(b),1.4(g), Proposition 3.5A , formula (18) and (19), we get
$$\begin{aligned}
&C_{r_0r_2r_3x_{i,j}u}C_{u^{-1}x_{k,l}r_3r_2r_0}\\
=&(C_{r_0}C_{r_2r_3x_{i,j}u}-C_{r_3x_{i,j}u}+\Box)(C_{u^{-1}x_{k,l}r_3r_2}C_{r_0}-C_{u^{-1}x_{k,l}r_3}+\Box)\\
=&C_{r_0}C_{r_2r_3x_{i,j}u}C_{u^{-1}x_{k,l}r_3r_2}C_{r_0}-C_{r_0}C_{r_2r_3x_{i,j}u}C_{u^{-1}x_{k,l}r_3}-C_{r_3x_{i,j}u}C_{u^{-1}x_{k,l}r_3r_2}C_{r_0}+C_{r_3x_{i,j}u}C_{u^{-1}x_{k,l}r_3}+\Box\\
=&\xi^6(\sum_{\gamma_{r_2r_3x_{i,j}u, u^{-1}x_{k,l}r_3r_2, ^{\star}z^{\star}}\ne 0} \gamma_{r_2r_3x_{i,j}u, u^{-1}x_{k,l}r_3r_2, ^{\star}z^{\star}}(C_{r_0}C_{^{\star}z^{\star}}C_{r_0})\\
&-\sum_{\gamma_{r_2r_3x_{i,j}u, u^{-1}x_{k,l}r_3, ^{\star}z^{\star*}}\ne 0}\gamma_{r_2r_3x_{i,j}u, u^{-1}x_{k,l}r_3, ^{\star}z^{\star*}}C_{r_0}C_{^{\star}z^{\star*}}\\
&-\sum_{\gamma_{r_3x_{i,j}u, u^{-1}x_{k,l}r_3r_2, ^{*\star}z^{\star}}\ne 0} \gamma_{r_3x_{i,j}u, u^{-1}x_{k,l}r_3r_2, ^{*\star}z^{\star}}C_{^{*\star}z^{\star}}C_{r_0}\\
&+\sum_{\gamma_{r_3x_{i,j}u, u^{-1}x_{k,l}r_3, ^{*\star}z^{\star*}}\ne 0} \gamma_{r_3x_{i,j}u, u^{-1}x_{k,l}r_3, ^{*\star}z^{\star*}}C_{^{*\star}z^{\star*}})+\triangle+\Box\\
=&\xi^6(\sum_{\gamma_{x_{i,j}, x_{k,l}, z}\ne 0}\gamma_{x_{i,j}, x_{k,l}, z}C_{r_0}C_{^{\star}z^{\star}}C_{r_0}-\sum_{\gamma_{x_{i,j}, x_{k,l}, z}\ne 0}\gamma_{x_{i,j}, x_{k,l}, z}}C_{r_0}C_{^{\star}z^{\star*}\\
&-\sum_{\gamma_{x_{i,j}, x_{k,l}, z}\ne 0}\gamma_{x_{i,j}, x_{k,l}, z}C_{^{*\star}z^{\star}}C_{r_0}+\sum_{\gamma_{x_{i,j}, x_{k,l}, z}\ne 0} \gamma_{x_{i,j}, x_{k,l}, z}C_{^{*\star}z^{\star*}}))+\triangle+\Box\end{aligned}$$
$$\begin{aligned}
=&\xi^6(\sum_{\gamma_{x_{i,j}, x_{k,l}, x_{m,n}}\ne 0}\gamma_{x_{i,j}, x_{k,l}, x_{m,n}}(C_{r_0r_2r_3x_{m,n}r_3r_2r_0}+C_{r_0r_2r_3x_{m,n}r_3}+C_{r_3x_{m,n}r_3r_2}C_{r_0})\\
&
-\sum_{\gamma_{x_{i,j}, x_{k,l}, x_{m,n}}\ne 0}\gamma_{x_{i,j}, x_{k,l}, x_{m,n}}(C_{r_0r_2r_3x_{m,n}r_3}+C_{r_3x_{m,n}r_3})-\sum_{\gamma_{x_{i,j}, x_{k,l}, x_{m,n}}\ne 0}\gamma_{x_{i,j}, x_{k,l}, x_{m,n}}C_{r_3x_{m,n}r_3r_2}C_{r_0}\\
&+\sum_{\gamma_{x_{i,j}, x_{k,l}, x_{m,n}}\ne 0}\gamma_{x_{i,j}, x_{k,l}, x_{m,n}}C_{r_3x_{m,n}r_3})+\triangle+\Box,\\
=&\xi^6\sum_{\gamma_{x_{i,j}, x_{k,l}, x_{m,n}}\ne 0}\gamma_{x_{i,j}, x_{k,l}, x_{m,n}}C_{r_0r_2r_3x_{m,n}r_3r_2r_0}+\triangle+\Box.
\end{aligned}$$

where $\triangle\in H$ is a sum of elements in the two-sided cell $c$ whose coefficient is of degree less than 6 and note $\gamma_{x_{i,j}, x_{k,l}, z}\ne 0$ only if $z=x_{m,n}$ for some $m, n\ge 0$.

Thus we have \begin{equation}
\gamma_{r_0\bar x_{i,j}, \bar x_{k,l}r_0, r_0(^{\star}z^{\star})r_0}=\gamma_{x_{i,j}, x_{k,l}, z}.\end{equation}

By 1.4(d) and the fact that any $x$ in $ \Gamma_{012}\cap \Gamma_{012}^{-1}$, we have $x=x^{-1}$, we get
\begin{align}\gamma_{\bar x_{k,l}r_0, r_0(^{\star}z^{\star})r_0, \bar x_{i,j}r_0}=\gamma_{x_{i,j}, x_{k,l}, z}.\end{align}

{\bf Step 2:} Compute $C_{r_0\bar x_{k,l}r_0}C_{r_0r_2r_3zr_3r_2r_0}$.

By formula (18), we get
\begin{equation}\begin{aligned}
C_{r_0\bar x_{k,l}r_0}C_{r_0(^{\star}z^{\star})r_0}=
C_{r_0}C_{\bar x_{k,l}r_0}C_{r_0(^{\star}z^{\star})r_0}-C_{r_3x_{k,l}r_3r_2r_0}C_{r_0(^{\star}z^{\star})r_0}+\Box.
\end{aligned}\end{equation}

Since by 1.4(d),1.4(g) and (20) we have 

$\gamma_{r_3x_{k,l}r_3r_2r_0, r_0(^{\star}z^{\star})r_0, r_3x_{i,j}r_3r_2r_0}=\gamma_{r_0r_2r_3x_{i,j}r_3, r_3x_{k,l}r_3r_2r_0, r_0(^{\star}z^{\star})r_0}=\gamma_{(r_0r_2r_3x_{i,j}r_3)^*, {}^*(r_3x_{k,l}r_3r_2r_0), r_0(^{\star}z^{\star})r_0}=\gamma_{r_0\bar x_{i,j}, \bar x_{k,l}r_0, r_0(^{\star}z^{\star})r_0}=\gamma_{x_{i,j}, x_{k,l}, z}.$

Similarly we get  $\gamma_{\bar x_{k,l}r_0, r_0(^{\star}z^{\star})r_0, \bar x_{i,j}r_0}=\gamma_{r_0\bar x_{i,j}, \bar x_{k,l}r_0, r_0(^{\star}z^{\star})r_0}=\gamma_{x_{i,j}, x_{k,l}, z}$.

Thus (22) turns into
$$\begin{aligned}
C_{r_0\bar x_{k,l}r_0}C_{r_0r_2r_3zr_3r_2r_0}&=C_{r_0}(C_{\bar x_{k,l}r_0}C_{r_0(^{\star}z^{\star})r_0})-C_{r_3x_{k,l}r_3r_2r_0}C_{r_0(^{\star}z^{\star})r_0}+\Box\\
&=\xi^6\sum_{\gamma_{x_{i,j},x_{k,l},z}\ne 0}\gamma_{x_{i,j},x_{k,l},z}(C_{r_0}C_{\bar x_{i,j}r_0}-C_{r_3x_{i,j}r_3r_2r_0})+\triangle+\Box\\
&=\xi^6\sum_{\gamma_{x_{i,j},x_{k,l},z}\ne 0}\gamma_{x_{i,j},x_{k,l},z}C_{r_0\bar x_{i,j}r_0}+\triangle+\Box,\end{aligned}$$
where $\triangle\in H$ is the sum of elements in $\Gamma_{02}$ whose coefficient is of degree less than 6.

Therefore, we get formula (18) $\gamma_{r_0\bar x_{i,j}r_0, r_0\bar x_{k,l}r_0, r_0r_2r_3zr_3r_2r_0}=\gamma_{x_{i,j}, x_{k,l}, z}.$ The proof is completed.
\end{proof}

{\bf Corollary 3.6B.} For any left cell $\Gamma,\ \Theta$ in $Y_4$,  we have a bijection
  $$\pi_{\Gamma \cap\Theta^{-1}}: \Gamma \cap\Theta^{-1}\to\text{Irr}(Sp_4(\mathbb C)\times\mathbb Z/2\mathbb Z),$$

 inducing an isomorphism of $\mathbb Z$-modules
  $$\pi_{\Gamma\cap\Theta^{-1}}: J_{\Gamma \cap\Theta^{-1}}\to\text{Rep}\ (Sp_4(\mathbb C)\times\mathbb Z/2\mathbb Z),$$
such that take any $x\in\Gamma \cap\Theta^{-1}$, $\Phi\in Y_4$, we have

(a) $\pi_{\Phi\cap\Theta^{-1}}(t_xt_y)=\pi_{\Gamma\cap\Theta^{-1}}(t_x)\pi_{\Phi\cap\Gamma^{-1}}(t_y), \text{for\ any}\ y\in \Phi\cap\Gamma^{-1}$,

(b) $\pi_{\Gamma\cap\Phi^{-1}}(t_zt_x)=\pi_{\Theta\cap\Phi^{-1}}(t_z)\pi_{\Gamma\cap\Theta^{-1}}(t_x), \text{for\ any}\ z\in \Theta\cap\Phi^{-1}$.

\begin{proof}
It follows Lemma 3.2 and subsection 3.3.
\end{proof}

  \medskip

  \noindent{\bf 3.7.}
In this subsection we establish a bijection

  $$\pi_{\Gamma_{03}\cap\Gamma_{03}^{-1}}: \Gamma_{03}\cap\Gamma_{03}^{-1}\to\text{Irr}(Sp_4(\mathbb C)\times\mathbb Z/2\mathbb Z)$$
  which induces an isomorphism of based rings
  $$\pi_{\Gamma_{03}\cap\Gamma_{03}^{-1}}: J_{\Gamma_{03}\cap \Gamma^{-1}_{03}}\to \text{Rep}(Sp_4(\mathbb C)\times\mathbb Z/2\mathbb Z).$$

  \medskip

Recall $x_{i,j}=w_{012}(\tau r_3r_2r_0r_1r_2)^i(r_3r_2r_0r_1)^{2j}$ for any $i, j\ge 0$. According to [D, Theorem 6.4] we have
$$ \begin{aligned}&\Gamma_{03}\cap \Gamma_{03}^{-1}\\
=&\{ r_3r_0r_2r_3x_{i,j}r_3r_2r_0r_3, \tau r_3r_1r_2r_3x_{i,j}r_3r_2r_0r_3 \,|\, i, j\ge 0\}\\
=&\{r_3xr_3\,|\, x\in\Gamma_{02}\cap\Gamma_{02}^{-1}\}.\end{aligned}$$

{\bf Proposition 3.7A:} The bijection
\begin{alignat*}{2} \pi_{\Gamma_{03}\cap\Gamma_{03}^{-1}}:    \Gamma_{03}\cap  \Gamma_{03}^{-1}&\longrightarrow \text{Irr}(Sp_4(\mathbb C)\times\mathbb Z/2\mathbb Z),\\
r_3xr_3&\longmapsto \pi_{\Gamma_{02}\cap\Gamma_{02}^{-1}}(x),
\end{alignat*}
   induces a based ring isomorphism
    $$\pi_{\Gamma_{03}\cap\Gamma_{03}^{-1}}: J_{\Gamma_{03}\cap \Gamma_{03}^{-1}}\longrightarrow \text{Rep}(Sp_4(\mathbb C)\times\mathbb Z/2\mathbb Z),\quad t_x\longmapsto \pi_{\Gamma_{03}\cap\Gamma_{03}^{-1}}(x).$$

Now we prove Proposition 3.7A.

    Let
    $$\hat x_{i,j}=r_0r_2r_3x_{i,j}r_3r_2r_0\in\Gamma_{02}\cap\Gamma_{02}^{-1}, \forall i, j\geq 0.$$

    Then any element in $\Gamma_{03}\cap\Gamma_{03}^{-1}$ is of the form $r_3\hat x_{i,j}r_3$ for some $i,j\ge 0$.

    Thanks to identity (c) in subsection 3.4 and the discussion in subsection 3.5 and subsection 3.6, to see Proposition 3.7A we only need to verify the following identities: for any $i,j,k,l\ge 0$, $z\in W$,
$$\gamma_{r_3\hat x_{i,j}r_3, r_3\hat x_{k,l}r_3, r_3zr_3}=\gamma_{\hat x_{i,j}, \hat x_{k,l}, z}.$$

\medskip

\begin{proof}
{\bf Step 1:} Compute $C_{r_3}C_{\hat x_{i,j}v}$, for any $v\in\{e, r_3\}$.

By 1.1(b),  we have $C_{r_3}C_{\hat x_{i,j}v}=C_{r_3\hat x_{i,j}v}+\sum\limits_{\substack{y\prec \hat x_{i,j}v\\ r_3y<y}}\mu(y, \hat x_{i,j}v)C_y.$
 Note that ${L}(\hat x_{i,j}v)=\{r_0,r_2\}$.

Assume $y\prec \hat x_{i,j}v$ and $r_3y<y$. Then
\begin{itemize}
\item If $r_0y>y$, then by 1.1(a) we get $r_0y=\hat x_{i,j}v$. This contradicts the assumption. So $r_0y>y$ would not occur.
\item If $r_2y>y$, then by 1.1(a) we get $r_2y=\hat x_{i,j}v$. This contradicts the assumption. So $r_2y>y$ would not occur.
\item If $r_0y<y, r_2y<y$ and $r_3y<y$, then $a(y)=a(w_{023})=9$ and $C_y\in H^{<w_{012}}$.
\end{itemize}

Therefore, we get
\begin{align}C_{r_3}C_{\hat x_{i,j}v}=C_{r_3\hat x_{i,j}v}+\Box.\end{align}

Hence
$$\begin{aligned}
C_{r_3\hat x_{i,j}}C_{\hat x_{k,l}r_3}=&C_{r_3}C_{\hat x_{i,j}}C_{\hat x_{k,l}}C_{r_3}+\Box\\
=&\xi^6\sum_{\gamma_{\hat x_{i,j},\hat x_{k,l},z}\ne 0}\gamma_{\hat x_{i,j},\hat x_{k,l},z}C_{r_3}C_zC_{r_3}+\triangle+\Box\\
=&\xi^6\sum_{\gamma_{\hat x_{i,j},\hat x_{k,l},z}\ne 0}\gamma_{\hat x_{i,j},\hat x_{k,l},z}C_{r_3zr_3}+\triangle+\Box,\end{aligned}$$
where $\triangle\in H$ is a sum of elements in $\Gamma_{03}\cap\Gamma_{03}^{-1}$ whose coefficient is of degree less than 6.

In other words, we get
\begin{align}\gamma_{r_3\hat x_{i,j}, \hat x_{k,l}r_3, r_3zr_3}=\gamma_{\hat x_{i,j}, \hat x_{k,l}, z}.\end{align}

By 1.4(d), we have \begin{align}\gamma_{\hat x_{k,l}r_3, r_3zr_3, \hat x_{i,j}r_3}=\gamma_{r_3\hat x_{i,j}, \hat x_{k,l}r_3, r_3zr_3}.\end{align}

\medskip

{\bf Step 2:} Compute $\gamma_{r_3\hat x_{k,l}r_3, r_3zr_3,r_3 \hat x_{i,j}r_3}$.

Note that $r_3\hat x_{i,j}r_3$ is the second element of the left string with respect to $\{r_2,r_3\}$, while $\hat x_{i,j}r_3$ is the first element of  the left string.

By 1.5(a), we get $\gamma_{r_3\hat x_{k,l}r_3, r_3zr_3,r_3 \hat x_{i,j}r_3}=\gamma_{\hat x_{k,l}r_3, r_3zr_3,\hat x_{i,j}r_3}+\gamma_{\hat x_{k,l}r_3, r_3zr_3,r_2r_3 \hat x_{i,j}r_3}.$

Since $r_2r_3 \hat x_{i,j}r_3\in(\Gamma'_{02})^{-1}$, by 1.4(c) we get $\gamma_{\hat x_{k,l}r_3, r_3zr_3,r_2r_3 \hat x_{i,j}r_3}=0$. Thus
\begin{align}\gamma_{r_3\hat x_{k,l}r_3, r_3zr_3,r_3 \hat x_{i,j}r_3}=\gamma_{\hat x_{k,l}r_3, r_3zr_3,\hat x_{i,j}r_3}.\end{align}

Combining formulas (24)-(26) and 1.4(d), the proof is completed.
\end{proof}

 \medskip

    {\bf Corollary 3.7B.} For any left cell $\Gamma,\ \Theta$ in $Y_2$, we have a bijection
  $$\pi_{\Gamma \cap\Theta^{-1}}: \Gamma \cap\Theta^{-1}\to\text{Irr}(Sp_4(\mathbb C)\times\mathbb Z/2\mathbb Z),$$
 inducing an isomorphism of $\mathbb Z$-modules
  $$\pi_{\Gamma\cap\Theta^{-1}}: J_{\Gamma \cap\Theta^{-1}}\to\text{Rep}\ (Sp_4(\mathbb C)\times\mathbb Z/2\mathbb Z),$$
such that take any $x\in\Gamma \cap\Theta^{-1}$, $\Phi\in Y_2$, we have

(a) $\pi_{\Phi\cap\Theta^{-1}}(t_xt_y)=\pi_{\Gamma\cap\Theta^{-1}}(t_x)\pi_{\Phi\cap\Gamma^{-1}}(t_y), \text{for\ any}\ y\in \Phi\cap\Gamma^{-1}$,

(b) $\pi_{\Gamma\cap\Phi^{-1}}(t_zt_x)=\pi_{\Theta\cap\Phi^{-1}}(t_z)\pi_{\Gamma\cap\Theta^{-1}}(t_x), \text{for\ any}\ z\in \Theta\cap\Phi^{-1}$.

  \begin{proof}
  It follows Lemma 3.2 and subsection 3.3.
  \end{proof}

  \medskip

\noindent{\bf 3.8.} In this subsection we
 establish a bijection
  $$ \pi_{\Gamma'_{013}\cap\Gamma_{013}^{'-1}}:  \Gamma'_{013}\cap\Gamma_{013}^{'-1}\to\text{Irr}(Sp_4(\mathbb C)\times\mathbb Z/2\mathbb Z)$$
  which induces an isomorphism of based rings
  $$ \pi_{\Gamma'_{013}\cap\Gamma_{013}^{'-1}}: J_{\Gamma'_{013}\cap \Gamma^{'-1}_{013}}\to \text{Rep}(Sp_4(\mathbb C)\times\mathbb Z/2\mathbb Z).$$

\medskip

  According to [D, Theorem 6.4] we have
$$ \Gamma'_{013}\cap \Gamma_{013}^{'-1}=\{ r_1r_3r_0r_2r_3w_{012}(r_3r_2r_0r_1r_2)^i(r_3r_2r_0r_1)^{2j}r_3r_2r_0r_3r_1\,|\, i, j\ge 0\}=\{r_1xr_1\,|\, x\in\Gamma_{03}\cap\Gamma_{03}^{-1}\}.$$

{\bf Proposition 3.8A:} The bijection
\begin{alignat*}{2} \pi_{\Gamma'_{013}\cap\Gamma_{013}^{'-1}}:    \Gamma'_{013}\cap  \Gamma_{013}^{'-1}&\longrightarrow \text{Irr}(Sp_4(\mathbb C)\times\mathbb Z/2\mathbb Z),\\
r_1xr_1&\longmapsto   \pi_{\Gamma_{03}\cap\Gamma_{03}^{-1}}(x) ,
\end{alignat*}
   induces a based ring isomorphism
    $$ \pi_{\Gamma'_{013}\cap\Gamma_{013}^{'-1}}: J_{ \Gamma'_{013}\cap \Gamma_{013}^{'-1}}\longrightarrow \text{Rep}(Sp_4(\mathbb C)\times\mathbb Z/2\mathbb Z),\quad t_x\longmapsto \pi_{\Gamma'_{013}\cap\Gamma_{013}^{'-1}}(x).$$

Now we prove  Proposition 3.8A.

    Let
    $$\begin{aligned}&x_{i,j}=w_{012}(r_3r_2r_0r_1r_2)^i(r_3r_2r_0r_1)^{2j},\\
    &\tilde x_{i,j}=r_3r_0r_2r_3x_{i,j}r_3r_2r_0r_3, \forall i, j\geq 0.\end{aligned}$$

    Then any element in $\Gamma'_{013}\cap\Gamma_{013}^{'-1}$ is of the form $r_1\tilde x_{i,j}r_1$ for some $i,j\ge 0$.

    Thanks to identity (c) in subsection 3.4 and the discussion in subsection 3.7, to see Proposition 3.8A we only need to verify the following identities: for any $i,j,k,l\ge 0$, $z\in W$,
$$\gamma_{r_1\tilde x_{i,j}r_1, r_1\tilde x_{k,l}r_1, r_1zr_1}=\gamma_{\tilde x_{i,j}, \tilde x_{k,l}, z}.$$

\begin{proof}
Let $\dot x_{i,j}=r_3r_0r_2r_3x_{i,j}$, $v\in\{e, r_3, r_3r_2, r_3r_2r_0, r_3r_2r_0r_3, r_3r_2r_0r_3r_1\}$.

{\bf Step 1:} Compute $C_{r_1}C_{\dot x_{i,j}v}$.

By 1.1(b),  we have $C_{r_1}C_{\dot x_{i,j}v}=C_{r_1\dot x_{i,j}v}+\sum\limits_{\substack{y\prec \dot x_{i,j}v\\ r_1y<y}}\mu(y, \dot x_{i,j}v)C_y.$
 Note that ${L}(\dot x_{i,j}v)=\{r_0,r_3\}$.

Assume $y\prec \dot x_{i,j}v$ and $r_1y<y$. Then
\begin{itemize}
\item If $r_0y>y$, then by 1.1(a) we get $r_0y=\dot  x_{i,j}v$. This contradicts the assumption. So $r_0y>y$ would not occur.
\item If $r_3y>y$, then by 1.1(a) we get $r_3y=\dot x_{i,j}v$. This contradicts the assumption. So $r_3y>y$ would not occur.
\item If $r_0y<y, r_1y<y$ and $r_3y<y$, then $y=r_0r_1r_3y'$ for some $y'\in W$ with $l(y)=l(y')+3$. If $r_2\in{L}(y)$, then $a(y)> a(w_{012})=6$ and $C_y\in H^{<w_{012}}$. Next assume $r_2\not\in{L}(y)$.
Note that $y$ is the first element of the left string with respect to $\{r_2,r_3\}$, while $\dot x_{i,j}v$ is the second element.

Then by 1.3(e) we get $\mu(y, \dot x_{i,j}v)=\mu(r_2y, r_2\dot x_{i,j}v)$ with ${L}(r_2\dot x_{i,j}v)=\{r_0,r_2\}$. If $r_0\in{L}(r_2y)$, then $r_2\in{L}(r_1r_3y')$ and $a(y)\ge a(r_1r_3y')=a(w_{123})=9$, indicating $C_y\in H^{<w_{012}}$.
Otherwise by 1.1(a) we get $r_0r_2y=r_2\dot x_{i,j}v$, contradicting $r_1y<y$.
So $r_0y<y, r_1y<y$ and $r_3y<y$ only deduce $C_y\in H^{<w_{012}}$.
\end{itemize}

Therefore,
\begin{align}C_{r_1}C_{\dot x_{i,j}v}=C_{r_1\dot x_{i,j}v}+\Box,\end{align}

and by 1.4(e) we get
\begin{align}C_{v(\dot x_{i,j})^{-1}}C_{r_1}=C_{v(\dot x_{i,j})^{-1}r_1}+\Box.\end{align}

Note $\tilde x_{i,j}=\dot x_{i,j}v$ with $v=r_3r_2r_0r_3$.
Hence by formula (27)-(28) and Proposition 3.7A, we get 
$$\begin{aligned}
C_{r_1\tilde x_{i,j}}C_{\tilde x_{k,l}r_1}=&C_{r_1}C_{\tilde x_{i,j}}C_{\tilde x_{k,l}}C_{r_1}+\Box\\
=&\xi^6\sum_{\gamma_{\tilde x_{i,j},\tilde x_{k,l},z}\ne 0}\gamma_{\tilde x_{i,j},\tilde x_{k,l},z}C_{r_1}C_zC_{r_1}+\triangle+\Box\\
=&\xi^6\sum_{\gamma_{\tilde x_{i,j},\tilde x_{k,l},z}\ne 0}\gamma_{\tilde x_{i,j},\tilde x_{k,l},z}C_{r_1zr_1}+\triangle+\Box,\end{aligned}$$
where $\triangle\in H$ is a sum of elements in $\Gamma'_{013}\cap(\Gamma'_{013})^{-1}$ whose coefficient is of degree less than 6. Note $C_{r_1}\triangle$ would not increase the degree of the coefficient.

In other words, we get
\begin{align}\gamma_{r_1\tilde x_{i,j}, \tilde x_{k,l}r_1, r_1zr_1}=\gamma_{\tilde x_{i,j}, \tilde x_{k,l}, z}.\end{align}

By 1.4(d), we have \begin{align}\gamma_{\tilde x_{k,l}r_1, r_1zr_1, \tilde x_{i,j}r_1}=\gamma_{r_1\tilde x_{i,j}, \tilde x_{k,l}r_1, r_1z^{-1}r_1}.\end{align}

{\bf Step 2:} Compute $C_{r_1}C_{\tilde x_{k,l}r_1}C_{ r_1zr_1}$.

By formula (30), (27) and Proposition 3.7A we get
$$\begin{aligned}
C_{r_1}C_{\tilde x_{k,l}r_1}C_{ r_1zr_1}
=&\xi^6\sum_{\gamma_{\tilde x_{k,l}r_1, r_1zr_1,\tilde x_{i,j}r_1}\ne 0}\gamma_{\tilde x_{k,l}r_1, r_1zr_1,\tilde x_{i,j}r_1}C_{r_1}C_{\tilde x_{i,j}r_1}+\triangle+\Box\\
=&\xi^6\sum_{\gamma_{\tilde x_{k,l}r_1, r_1zr_1,\tilde x_{i,j}r_1}\ne 0}\gamma_{\tilde x_{k,l}r_1, r_1zr_1,\tilde x_{i,j}r_1}C_{r_1\tilde x_{i,j}r_1}+\triangle+\Box\\
=&C_{r_1 \tilde x_{k,l}r_1}C_{r_1zr_1}+\triangle+\Box\\
&(\text{by the composition law of product}),\end{aligned}$$

where $\triangle\in H$ is the sum of elements in $\Gamma'_{013}$ whose coefficient is of degree less than 6.

Hence \begin{align} \gamma_{r_1 \tilde x_{k,l}r_1, r_1zr_1, r_1\tilde x_{i,j}r_1}=\gamma_{\tilde x_{k,l}r_1, r_1zr_1,\tilde x_{i,j}r_1}.\end{align}

Combining formulas (29)-(31) and 1.4(d), the proof is completed.
\end{proof}

  {\bf Corollary 3.8B.} For any left cell $\Gamma,\ \Theta$ in $Y_3$, we have a bijection
  $$\pi_{\Gamma\cap\Theta^{-1}}: \Gamma \cap\Theta^{-1}\to\text{Irr}(Sp_4(\mathbb C)\times\mathbb Z/2\mathbb Z),$$
inducing an isomorphism of $\mathbb Z$-modules
  $$\pi_{\Gamma\cap\Theta^{-1}}: J_{\Gamma \cap\Theta^{-1}}\to\text{Rep}\ (Sp_4(\mathbb C)\times\mathbb Z/2\mathbb Z),$$
such that take any $x\in\Gamma \cap\Theta^{-1}$, $\Phi\in Y_3$, we have

(a) $\pi_{\Phi\cap\Theta^{-1}}(t_xt_y)=\pi_{\Gamma\cap\Theta^{-1}}(t_x)\pi_{\Phi\cap\Gamma^{-1}}(t_y), \text{for\ any}\ y\in \Phi\cap\Gamma^{-1}$,

(b) $\pi_{\Gamma\cap\Phi^{-1}}(t_zt_x)=\pi_{\Theta\cap\Phi^{-1}}(t_z)\pi_{\Gamma\cap\Theta^{-1}}(t_x), \text{for\ any}\ z\in \Theta\cap\Phi^{-1}$.

  \begin{proof}
  It follows Lemma 3.2 and subsection 3.3.
  \end{proof}

Combining Proposition 3.5A, Proposition 3.6A, Proposition 3.7A and Proposition 3.8A, we have proved Theorem 3.4A.

Next we prove Theorem 3.4B.

  \medskip
  
\noindent{\bf 3.9.}  In this subsection we establish six bijections
  $$\pi_{\Gamma_{012}\cap\Gamma_{02}^{-1}}: \Gamma_{012}\cap\Gamma_{02}^{-1}\to\text{the set of isomorphism classes of irreducible $F_c$-v.b. on}\ \Gamma_{02}\times \Gamma_{012};$$
   $$\pi_{ \Gamma_{012}\cap\Gamma_{03}^{-1}}: \Gamma_{012}\cap\Gamma_{03}^{-1}\to\text{the set of isomorphism classes of irreducible $F_c$-v.b. on}\ \Gamma_{03}\times \Gamma_{012};$$
 $$\pi_{\Gamma_{012}\cap\Gamma_{013}^{'-1}}: \Gamma_{012}\cap\Gamma_{013}^{'-1}\to\text{the set of isomorphism classes of irreducible $F_c$-v.b. on}\ \Gamma'_{013}\times \Gamma_{012} ;$$
$$\pi_{\Gamma_{02}\cap\Gamma_{03}^{-1}}: \Gamma_{02}\cap\Gamma_{03}^{-1}\to\text{the set of isomorphism classes of irreducible $F_c$-v.b. on}\ \Gamma_{03}\times\Gamma_{02};$$
$$\pi_{\Gamma_{02}\cap\Gamma_{013}^{'-1}}: \Gamma_{02}\cap\Gamma_{013}^{'-1}\to\text{the set of isomorphism classes of irreducible $F_c$-v.b. on}\ \Gamma'_{013}\times \Gamma_{02};$$
$$\pi_{ \Gamma_{03}\cap\Gamma_{013}^{'-1}}: \Gamma_{03}\cap\Gamma_{013}^{'-1}\to\text{the set of isomorphism classes of irreducible $F_c$-v.b. on}\ \Gamma'_{013}\times \Gamma_{03}.$$

\medskip

  The inverses give bijections
  $$\pi_{\Gamma_{02}\cap\Gamma_{012}^{-1}}: \Gamma_{02}\cap\Gamma_{012}^{-1}\to\text{the set of isomorphism classes of irreducible $F_c$-v.b. on}\ \Gamma_{012}\times \Gamma_{02};$$
    $$\pi_{\Gamma_{03}\cap\Gamma_{012}^{-1}}: \Gamma_{03}\cap\Gamma_{012}^{-1}\to\text{the set of isomorphism classes of irreducible $F_c$-v.b. on}\ \Gamma_{012}\times \Gamma_{03};$$
  $$\pi_{\Gamma'_{013}\cap\Gamma_{012}^{-1}}: \Gamma'_{013}\cap\Gamma_{012}^{-1}\to\text{the set of isomorphism classes of irreducible $F_c$-v.b. on}\ \Gamma_{012}\times\Gamma'_{013};$$
    $$\pi_{\Gamma_{03}\cap\Gamma_{02}^{-1}}: \Gamma_{03}\cap\Gamma_{02}^{-1}\to\text{the set of isomorphism classes of irreducible $F_c$-v.b. on}\ \Gamma_{02} \times \Gamma_{03};$$
    $$\pi_{\Gamma'_{013}\cap\Gamma_{02}^{-1}}: \Gamma'_{013}\cap\Gamma_{02}^{-1}\to\text{the set of isomorphism classes of irreducible $F_c$-v.b. on}\ \Gamma_{02}\times \Gamma'_{013};$$
   $$\pi_{\Gamma'_{013}\cap\Gamma_{03}^{-1}}: \Gamma'_{013}\cap\Gamma_{03}^{-1}\to\text{the set of isomorphism classes of irreducible $F_c$-v.b. on}\ \Gamma_{03}\times \Gamma'_{013}.$$
  
  According to [D, Theorem 6.4] and the notation in 3.5, we have
$$\begin{aligned} \Gamma_{012}\cap \Gamma_{02}^{-1}&=\{ r_0r_2r_3w_{012}(r_3r_2r_0r_1r_2)^i(r_3r_2r_0r_1)^{2j}\,|\, i, j\ge 0\}=\{r_0({}^\star x)\,|\, x\in\Gamma_{012}\cap\Gamma_{012}^{-1}\};\\
\Gamma_{012}\cap \Gamma_{03}^{-1}&=\{ r_3r_0r_2r_3w_{012}(r_3r_2r_0r_1r_2)^i(r_3r_2r_0r_1)^{2j}\,|\, i, j\ge 0\}
=\{r_3x\, |\, x\in\Gamma_{012}\cap\Gamma_{02}^{-1}\};\\
\Gamma_{012}\cap \Gamma_{013}^{'-1}&=\{ r_1r_3r_0r_2r_3w_{012}(r_3r_2r_0r_1r_2)^i(r_3r_2r_0r_1)^{2j}\,|\, i, j\ge 0\}=\{r_1x\, |\, x\in\Gamma_{012}\cap\Gamma_{03}^{-1}\};\\
\Gamma_{02}\cap \Gamma_{03}^{-1}&=\{ r_3r_0r_2r_3w_{012}(r_3r_2r_0r_1r_2)^i(r_3r_2r_0r_1)^{2j}r_3r_2r_0\,|\, i, j\ge 0\}=\{r_3x\, |\, x\in\Gamma_{02}\cap\Gamma_{02}^{-1}\};\\
\Gamma_{02}\cap\Gamma_{013}^{'-1}&=\{ r_1r_3r_0r_2r_3w_{012}(r_3r_2r_0r_1r_2)^i(r_3r_2r_0r_1)^{2j}r_3r_2r_0\,|\, i, j\ge 0\}=\{r_1x\, |\, x\in\Gamma_{02}\cap\Gamma_{03}^{-1}\};\\
\Gamma_{03}\cap \Gamma_{013}^{'-1}&=\{ r_1r_3r_0r_2r_3w_{012}(r_3r_2r_0r_1r_2)^i(r_3r_2r_0r_1)^{2j}r_3r_2r_0r_3\,|\, i, j\ge 0\}=\{r_1x\, |\, x\in\Gamma_{03}\cap\Gamma_{03}^{-1}\}.\end{aligned}$$

{\bf Proposition 3.9A:} 
\begin{itemize}
\item[(i)] The bijection
\begin{alignat*}{2} \pi_{ \Gamma_{012}\cap  \Gamma_{02}^{-1}}:    \Gamma_{012}\cap  \Gamma_{02}^{-1}&\longrightarrow  \text{Irr}(Sp_4(\mathbb C)\times\mathbb Z/2\mathbb Z),\\
r_0(^\star x)&\longmapsto \pi_{\Gamma_{012}\cap  \Gamma_{012}^{-1}}(x)
\end{alignat*}
induces an isomorphism of $\mathbb Z$-modules
$$\pi_{\Gamma_{012}\cap \Gamma^{-1}_{02}}: J_{\Gamma_{012}\cap \Gamma^{-1}_{02}}\to  K_{F}( \Gamma_{02}\times\Gamma_{012})\cong\text{Rep}\ (Sp_4(\mathbb C)\times\mathbb Z/2\mathbb Z) ), t_x\mapsto \pi_{\Gamma_{012}\cap \Gamma^{-1}_{02}}(x),$$

\item[(ii)]the bijection
\begin{alignat*}{2} \pi_{\Gamma_{012}\cap  \Gamma_{03}^{-1}}:    \Gamma_{012}\cap  \Gamma_{03}^{-1}&\longrightarrow  \text{Irr}(Sp_4(\mathbb C)\times\mathbb Z/2\mathbb Z),\\
r_3x&\longmapsto \pi_{\Gamma_{012}\cap\Gamma_{02}^{-1}}(x)
\end{alignat*}
induces an isomorphism of $\mathbb Z$-modules
$$\pi_{\Gamma_{012}\cap \Gamma^{-1}_{03}}: J_{\Gamma_{012}\cap \Gamma^{-1}_{03}}\to  K_{F}( \Gamma_{03}\times\Gamma_{012})\cong\text{Rep}\ (Sp_4(\mathbb C)\times\mathbb Z/2\mathbb Z) ), t_x\mapsto \pi_{\Gamma_{012}\cap \Gamma^{-1}_{03}}(x),$$

\item[(iii)] the bijection
\begin{alignat*}{2} \pi_ {\Gamma_{012}\cap  (\Gamma'_{013})^{-1}}:    \Gamma_{012}\cap  (\Gamma'_{013})^{-1}&\longrightarrow  \text{Irr}(Sp_4(\mathbb C)\times\mathbb Z/2\mathbb Z),\\
r_1x&\longmapsto \pi_ {\Gamma_{012}\cap  (\Gamma_{03})^{-1}}(x).
\end{alignat*}
induces an isomorphism of $\mathbb Z$-modules
$$\pi_{\Gamma_{012}\cap \Gamma_{013}^{'-1}}: J_{\Gamma_{012}\cap \Gamma_{013}^{'-1}}\to  K_{F}( \Gamma'_{013}\times\Gamma_{012})\cong\text{Rep}\ (Sp_4(\mathbb C)\times\mathbb Z/2\mathbb Z) ), t_x\mapsto \pi_{\Gamma_{012}\cap\Gamma_{013}^{'-1}}(x),$$

\item[(iv)]the bijection
\begin{alignat*}{2} \pi_{\Gamma_{02}\cap  \Gamma_{03}^{-1}}:    \Gamma_{02}\cap  \Gamma_{03}^{-1}&\longrightarrow  \text{Irr}(Sp_4(\mathbb C)\times\mathbb Z/2\mathbb Z),\\
r_3x&\longmapsto \pi_{\Gamma_{02}\cap  \Gamma_{02}^{-1}}(x).
\end{alignat*}
induces an isomorphism of $\mathbb Z$-modules
$$\pi_{\Gamma_{02}\cap \Gamma^{-1}_{03}}: J_{\Gamma_{02}\cap \Gamma^{-1}_{03}}\to  K_{F}( \Gamma_{03}\times\Gamma_{02})\cong\text{Rep}\ (Sp_4(\mathbb C)\times\mathbb Z/2\mathbb Z) ), t_x\mapsto \pi_{\Gamma_{02}\cap\Gamma_{03}^{-1}}(x),$$

\item[(v)]the bijection
\begin{alignat*}{2} \pi_{\Gamma_{02}\cap  (\Gamma'_{013})^{-1}}:    \Gamma_{02}\cap  (\Gamma'_{013})^{-1}&\longrightarrow  \text{Irr}(Sp_4(\mathbb C)\times\mathbb Z/2\mathbb Z),\\
r_1x&\longmapsto   \pi_{\Gamma_{02}\cap \Gamma_{03}^{-1}}(x),
\end{alignat*}
induces an isomorphism of $\mathbb Z$-modules
$$\pi_{\Gamma_{02}\cap \Gamma^{'-1}_{013}}: J_{\Gamma_{02}\cap \Gamma^{'-1}_{013}}\to  K_{F}( \Gamma'_{013}\times\Gamma_{02})\cong\text{Rep}\ (Sp_4(\mathbb C)\times\mathbb Z/2\mathbb Z) ), t_x\mapsto \pi_{\Gamma_{02}\cap\Gamma_{013}^{'-1}}(x),$$

\item[(vi)]the bijection
\begin{alignat*}{2} \pi_{\Gamma_{03}\cap  (\Gamma'_{013})^{-1}}:    \Gamma_{03}\cap  (\Gamma'_{013})^{-1}&\longrightarrow  \text{Irr}(Sp_4(\mathbb C)\times\mathbb Z/2\mathbb Z),\\
r_1x&\longmapsto \pi_{\Gamma_{03}\cap \Gamma_{03}^{-1}}(x).
\end{alignat*}
induces an isomorphism of $\mathbb Z$-modules
$$\pi_{\Gamma_{03}\cap \Gamma^{'-1}_{013}}: J_{\Gamma_{02}\cap \Gamma^{'-1}_{013}}\to  K_{F}( \Gamma'_{013}\times\Gamma_{03})\cong\text{Rep}\ (Sp_4(\mathbb C)\times\mathbb Z/2\mathbb Z) ), t_x\mapsto \pi_{\Gamma_{03}\cap\Gamma_{013}^{'-1}}(x),$$
\end{itemize}
such that for any $u\in\Gamma_{012}\cap \Gamma^{-1}_{02}$, $v\in\Gamma_{012}\cap\Gamma_{03}^{-1}$, $w\in\Gamma_{012}\cap\Gamma_{013}^{'-1}$, $x\in\Gamma_{02}\cap\Gamma_{03}^{-1}$, $y\in\Gamma_{02}\cap\Gamma_{013}^{'-1}$ and $z\in\Gamma_{03}\cap\Gamma_{013}^{'-1}$,
we have
\begin{itemize}
\item[(a)]$\pi_{\Gamma\cap \Gamma^{-1}_{02}}(t_ut_m)=\pi_{\Gamma_{012}\cap \Gamma^{-1}_{02}}(t_u)\pi_{\Gamma\cap\Gamma_{012}^{-1}}(t_m)$, for any
$m\in\Gamma\cap\Gamma_{012}^{-1}$, $\Gamma\in Y_1\cup Y_2\cup Y_3\cup Y_4$;
\item[(b)] $\pi_{\Gamma_{012}\cap \Theta^{-1}}(t_pt_u)=\pi_{\Gamma_{02}\cap\Theta^{-1}}(t_p)\pi_{\Gamma_{012}\cap \Gamma^{-1}_{02}}(t_u)$, for any $p\in\Gamma_{02}\cap\Theta^{-1}$, $\Theta\in Y_1\cup Y_2\cup Y_3\cup Y_4$;
\item[(c)] $\pi_{\Gamma\cap \Gamma^{-1}_{03}}(t_vt_m)=\pi_{\Gamma_{012}\cap \Gamma^{-1}_{03}}(t_v)\pi_{\Gamma\cap\Gamma_{012}^{-1}}(t_m)$, for any
$m\in\Gamma\cap\Gamma_{012}^{-1}$, $\Gamma\in Y_1\cup Y_2\cup Y_3\cup Y_4$;
\item[(d)] $\pi_{\Gamma_{012}\cap \Theta^{-1}}(t_pt_v)=\pi_{\Gamma_{03}\cap\Theta^{-1}}(t_p)\pi_{\Gamma_{012}\cap \Gamma^{-1}_{03}}(t_v)$, for any
$p\in\Gamma_{03}\cap\Theta^{-1}$, $\Theta\in Y_1\cup Y_2\cup Y_3\cup Y_4$;
\item[(e)] $\pi_{\Gamma\cap (\Gamma'_{013})^{-1}}(t_wt_m)=\pi_{\Gamma_{012}\cap (\Gamma'_{013})^{-1}}(t_w)\pi_{\Gamma_{012}\cap\Gamma_{012}^{-1}}(t_m)$, for any
$m\in\Gamma\cap\Gamma_{012}^{-1}$, $\Gamma\in Y_1\cup Y_2\cup Y_3\cup Y_4$;
\item[(f)] $\pi_{\Gamma_{012}\cap \Theta^{-1}}(t_pt_w)=\pi_{\Gamma'_{013}\cap\Theta^{-1}}(t_p)\pi_{\Gamma_{012}\cap (\Gamma'_{013})^{-1}}(t_w)$, for any
$p\in\Gamma'_{013}\cap\Theta^{-1}$, $\Theta\in Y_1\cup Y_2\cup Y_3\cup Y_4$;
\item[(g)] $\pi_{\Gamma\cap \Gamma_{03}^{-1}}(t_xt_m)=\pi_{\Gamma_{02}\cap\Gamma_{03}^{-1}}(t_x)\pi_{\Gamma\cap\Gamma_{02}^{-1}}(t_m)$, for any
$m\in\Gamma\cap\Gamma_{02}^{-1}$, $\Gamma\in Y_1\cup Y_2\cup Y_3\cup Y_4$;
\item[(h)] $\pi_{\Gamma_{02}\cap \Theta^{-1}}(t_pt_x)=\pi_{\Gamma_{03}\cap\Theta^{-1}}(t_p)\pi_{\Gamma_{02}\cap\Gamma_{03}^{-1}}(t_x)$, for any
$p\in\Gamma_{03}\cap\Theta^{-1}$, $\Theta\in Y_1\cup Y_2\cup Y_3\cup Y_4$;
\item[(i)] $\pi_{\Gamma\cap \Gamma_{013}^{'-1}}(t_yt_m)=\pi_{\Gamma_{02}\cap\Gamma_{013}^{'-1}}(t_y)\pi_{\Gamma\cap\Gamma_{02}^{-1}}(t_m)$, for any
$m\in\Gamma\cap\Gamma_{02}^{-1}$, $\Gamma\in Y_1\cup Y_2\cup Y_3\cup Y_4$;
\item[(j)] $\pi_{\Gamma_{02}\cap \Theta^{-1}}(t_pt_y)=\pi_{\Gamma'_{013}\cap\Theta^{-1}}(t_p)\pi_{\Gamma_{02}\cap\Gamma_{013}^{'-1}}(t_y)$, for any
$p\in\Gamma'_{013}\cap\Theta^{-1}$, $\Theta\in Y_1\cup Y_2\cup Y_3\cup Y_4$;
\item[(k)] $\pi_{\Gamma\cap \Gamma_{013}^{'-1}}(t_zt_m)=\pi_{\Gamma_{03}\cap\Gamma_{013}^{'-1}}(t_z)\pi_{\Gamma\cap\Gamma_{03}^{-1}}(t_m)$, for any
$m\in\Gamma\cap\Gamma_{03}^{-1}$, $\Gamma\in Y_1\cup Y_2\cup Y_3\cup Y_4$;
\item[(l)] $\pi_{\Gamma_{03}\cap \Theta^{-1}}(t_pt_z)=\pi_{\Gamma'_{013}\cap\Theta^{-1}}(t_p)\pi_{\Gamma_{03}\cap\Gamma_{013}^{'-1}}(t_z)$, for any
$p\in\Gamma'_{013}\cap\Theta^{-1}$, $\Theta\in Y_1\cup Y_2\cup Y_3\cup Y_4$.
\end{itemize}

\medskip

\begin{proof}
Note $r_2r_3x_{i,j}={}^\star x_{i,j}$,  to prove assertion (a) and (b) of Proposition 3.9A, we need to check the following identities:
\begin{align}
&\gamma_{r_0(^\star x_{i,j}), x_{k,l}, h_1}\ne 0\ \text{only\ if}\ h_1=r_0(^\star z), \ \text{and}\ \gamma_{r_0(^\star x_{i,j}), x_{k,l}, r_0(^\star z)}=\gamma_{x_{i,j}, x_{k,l}, z};\\
&\gamma_{r_0(^\star x_{i,j}), (x_{k,l}^\star)r_0, h_2}\ne 0\ \text{ only\ if}\ h_2=r_0(^\star z^\star)r_0, \ \text{and}\ \gamma_{r_0(^\star x_{i,j}), (x_{k,l}^\star)r_0, r_0(^\star z^\star)r_0}=\gamma_{x_{i,j}, x_{k,l}, z};\end{align}
\begin{align}
&\gamma_{r_0(^\star x_{i,j}), (x_{k,l}^\star)r_0r_3, h_3}\ne 0\ \text{only\ if}\ h_3=r_0(^\star z^\star)r_0r_3, \ \text{and}\  \gamma_{r_0(^\star x_{i,j}), (x_{k,l}^\star)r_0r_3, r_0(^\star z^\star)r_0r_3}=\gamma_{x_{i,j}, x_{k,l}, z};\\
&\gamma_{r_0(^\star x_{i,j}), (x_{k,l}^\star)r_0r_3r_1, h_4}\ne 0\ \text{only\ if}\ h_4=r_0(^\star z^\star)r_0r_3r_1, \ \text{and}\ \gamma_{r_0(^\star x_{i,j}), (x_{k,l}^\star)r_0r_3r_1, r_0(^\star z^\star)r_0r_3r_1}=\gamma_{x_{i,j}, x_{k,l}, z};\\
&\gamma_{(x_{i,j}^\star)r_0, r_0(^\star x_{k,l}), h_5}=\gamma_{x_{i,j}, x_{k,l}, h_5};\\
&\gamma_{r_0(^\star x_{i,j}^\star)r_0, r_0(^\star x_{k,l}), h_6}\ne 0\ \text{only\ if}\ h_6=r_0(^\star z), \ \text{and}\  \gamma_{r_0(^\star x_{i,j}^\star)r_0, r_0(^\star x_{k,l}), h_6}=\gamma_{x_{i,j}, x_{k,l}, z};\\
&\gamma_{r_3r_0(^\star x_{i,j}^\star)r_0, r_0(^\star x_{k,l}), h_7}\ne 0\ \text{only\ if}\ h_7=r_3r_0(^\star z), \ \text{and}\  \gamma_{r_3r_0(^\star x_{i,j}^\star)r_0, r_0(^\star x_{k,l}), h_7}=\gamma_{x_{i,j}, x_{k,l}, z};\\
&\gamma_{r_1r_3r_0(^\star x_{i,j}^\star)r_0, r_0(^\star x_{k,l}), h_8}\ne 0\ \text{only\ if}\ h_8=r_1r_3r_0(^\star z), \ \text{and}\  \gamma_{r_1r_3r_0(^\star x_{i,j}^\star)r_0, r_0(^\star x_{k,l}), h_8}=\gamma_{x_{i,j}, x_{k,l}, z};
\end{align}

\medskip

For (32), note $r_2r_3x_{i,j}={}^\star x_{i,j}$ and $r_3x_{i,j}={}^{*\star} x_{i,j}$.

Then by formula (18),1.4(g) and Proposition 3.5A , we get
$$\begin{aligned}&C_{r_0({}^\star x_{i,j})}C_{x_{k,l}}\\
=&C_{r_0r_2r_3x_{i,j}}C_{x_{k,l}}\\=&C_{r_0}C_{r_2r_3x_{i,j}}C_{x_{k,l}}-C_{r_3x_{i,j}}C_{x_{k,l}}+\Box\\
=&\xi^6(\sum_{\gamma_{r_2r_3x_{i,j}, x_{k,l},h}\ne 0}\gamma_{r_2r_3x_{i,j}, x_{k,l},h}C_{r_0}C_{h}-\sum _{\gamma_{r_3x_{i,j}, x_{k,l}, y}\ne 0}\gamma_{r_3x_{i,j}, x_{k,l}, y}C_{y})+\triangle+\Box\\
=&\xi^6(\sum_{\gamma_{^\star x_{i,j}, x_{k,l},h}\ne 0}\gamma_{^\star x_{i,j}, x_{k,l},h}C_{r_0}C_{h}-\sum _{\gamma_{^{*\star}x_{i,j}, x_{k,l}, y}\ne 0}\gamma_{^{*\star}x_{i,j}, x_{k,l}, y}C_{y})+\triangle+\Box\\
=&\xi^6(\sum_{\gamma_{x_{i,j}, x_{k,l},^\star h}\ne 0}\gamma_{x_{i,j}, x_{k,l},^\star h}C_{r_0}C_{h}-\sum _{\gamma_{x_{i,j}, x_{k,l}, ^{*\star}y}\ne 0}\gamma_{x_{i,j}, x_{k,l}, ^{*\star}y}C_{y})+\triangle+\Box\\
=&\xi^6(\sum_{\gamma_{x_{i,j}, x_{k,l}, x_{m,n}}\ne 0}\gamma_{x_{i,j}, x_{k,l},x_{m,n}}(C_{r_0s_2s_3x_{m,n}}+C_{r_3x_{m,n}})-\sum _{\gamma_{x_{i,j}, x_{k,l}, x_{m,n}}\ne 0}\gamma_{x_{i,j}, x_{k,l}, x_{m,n}}C_{r_3x_{m,n}})+\triangle+\Box\\
=&\xi^6\sum_{\gamma_{x_{i,j}, x_{k,l}, z}\ne 0}\gamma_{x_{i,j}, x_{k,l}, z}C_{r_0(^\star z)}+\triangle+\Box.\end{aligned}$$

where $\triangle\in H$ is the sum of elements in $B$ whose coefficient is of degree less than 6.

Then (32) is valid.

For (33), since $x_{k,l}r_2r_3=x_{k,l}^\star$, combining with (32), by 1.4(g) it suffices to show \begin{align}\gamma_{r_0(^\star x_{i,j}), (x_{k,l}^\star)r_0, r_0(^\star z^\star)r_0}=\gamma_{r_0(^\star x_{i,j}), x_{k,l}, r_0(^\star z)}.\end{align}

By (19) and 1.4(g), we get
$$\begin{aligned}C_{r_0(^\star x_{i,j})}C_{(x_{k,l}^\star) r_0}=&C_{r_0(^\star x_{i,j})}C_{x_{k,l}^\star}C_{r_0}-C_{r_0(^\star x_{i,j})}C_{x_{k,l}^{\star*}}+\Box\\
=&\xi^6(\sum_{\gamma_{r_0(^\star x_{i,j}), x_{k,l}^\star, h}\ne 0}\gamma_{r_0(^\star x_{i,j}), x_{k,l}^\star, h}C_{h}C_{r_0} -\sum_{\gamma_{r_0(^\star x_{i,j}), x_{k,l}^{\star*}, y}\ne 0}\gamma_{r_0(^\star x_{i,j}), x_{k,l}^{\star*}, y}C_{y})\\
&+\triangle+\Box\\
=&\xi^6(\sum_{\gamma_{r_0(^\star x_{i,j}), x_{k,l}, h^\star}\ne 0}\gamma_{r_0(^\star x_{i,j}), x_{k,l}, h^\star}C_{h}C_{r_0}-\sum_{\gamma_{r_0(^\star x_{i,j}), x_{k,l}, y^{*\star}}\ne 0}\gamma_{r_0(^\star x_{i,j}), x_{k,l}, y^{*\star}}C_{y})\\
&+\triangle+\Box\\
=&\xi^6(\sum_{\gamma_{r_0(^\star x_{i,j}), x_{k,l}, r_0(^\star z)}\ne 0}\gamma_{r_0(^\star x_{i,j}), x_{k,l}, r_0(^\star z)}C_{r_0(^\star z^\star)}C_{r_0}\\
&\ \ \ \ -\sum_{\gamma_{r_0(^\star x_{i,j}), x_{k,l}, r_0(^\star z)}\ne 0}\gamma_{r_0(^\star x_{i,j}), x_{k,l},r_0(^\star z)}C_{r_0(^\star z^{\star*})})+\triangle+\Box\\
=&\xi^6(\sum_{\gamma_{r_0(^\star x_{i,j}), x_{k,l}, r_0(^\star x_{m,n})}\ne 0}\gamma_{r_0(^\star x_{i,j}), x_{k,l}, r_0(^\star x_{m,n})}(C_{r_0(^\star x_{m,n}^
\star)r_0}-C_{r_0(^\star x_{m,n})s_3})\\
&\ \ \ \ -\sum_{\gamma_{r_0(^\star x_{i,j}), x_{k,l}, r_0(^\star x_{m,n})}\ne 0}\gamma_{r_0(^\star x_{i,j}), x_{k,l},r_0(^\star x_{m,n})}C_{r_0(^\star x_{m,n})s_3})+\triangle+\Box\\
=&\xi^6\sum_{\gamma_{r_0(^\star x_{i,j}), x_{k,l}, r_0(^\star x_{m,n})}\ne 0}\gamma_{r_0(^\star x_{i,j}), x_{k,l}, r_0(^\star x_{m,n})}C_{r_0(^\star x_{m,n}^
\star)r_0}+\triangle+\Box\\
=&\xi^6\sum_{\gamma_{r_0(^\star x_{i,j}), x_{k,l}, r_0(^\star z)}\ne 0}\gamma_{r_0(^\star x_{i,j}), x_{k,l}, r_0(^\star z)}C_{r_0(^\star z^
\star)r_0}+\triangle+\Box
\end{aligned}$$

where $\triangle\in H$ is the sum of elements in $B$ whose coefficient is of degree less than 6.

Then (40) is true and (33) is valid.

For (34), by (33) it suffices to show \begin{align}\gamma_{r_0(^\star x_{i,j}), (x_{k,l}^{\star})r_0r_3, r_0(^\star z^\star)r_0r_3}=\gamma_{r_0(^\star x_{i,j}), (x_{k,l}^{\star})r_0, r_0(^\star z^\star)r_0}.\end{align}

By 1.5(b), (33) and Proposition 3.5A, we get

$\gamma_{x_{i,j}, x_{k,l}, z}\ne 0\ \text{if\ and\ only\ if}\ z=x_{m,n}, $ and

$\begin{aligned}\gamma_{r_0(^\star x_{i,j}), x_{k,l}^{\star}r_0, r_0(^\star x_{m,n}^\star)r_0}&=\gamma_{r_0r_2r_3x_{i,j}, x_{k,l}r_3r_2r_0, r_0r_2r_3x_{m,n}r_3r_2r_0}\\
&=\gamma_{r_0r_2r_3x_{i,j}, x_{k,l}r_3r_2r_0r_3, r_0r_2r_3x_{m,n}r_3r_2r_0r_3}-\gamma_{r_0r_2r_3x_{i,j}, x_{k,l}r_3r_2r_0, r_0r_2r_3x_{m,n}r_3r_2r_0r_3r_2}\\
&=\gamma_{r_0(^\star x_{i,j}), (x_{k,l}^\star) r_0r_3, r_0(^\star x_{m,n}^\star)r_0r_3}-\gamma_{r_0(^\star x_{i,j}), x_{k,l}^\star r_0, r_0(^\star x_{m,n}^\star)r_0r_3r_2}\\
&=\gamma_{r_0(^\star x_{i,j}), (x_{k,l}^\star) r_0r_3, r_0(^\star x_{m,n}^\star)r_0r_3},\end{aligned}$
 so (41) is valid and (34) is true.

For (35), it suffices to show \begin{align}\gamma_{r_0(^\star x_{i,j}), (x_{k,l}^\star )r_0r_3r_1, r_0(^\star z^\star)r_0r_3r_1}=\gamma_{r_0(^\star x_{i,j}), (x_{k,l}^\star) r_0r_3,r_0(^\star z^\star)r_0r_3}.\end{align}

By (28), (34),1.4(b) and Proposition 3.5A, we get
\begin{equation}\begin{aligned}
&C_{r_0(^\star x_{i,j})}C_{(x_{k,l}^\star)r_0r_3r_1}
=C_{r_0r_2r_3x_{i,j}}C_{x_{k,l}r_3r_2r_0r_3r_1}
=C_{r_0r_2r_3x_{i,j}}C_{x_{k,l}r_3r_2r_0r_3}C_{r_1}+\Box\\
=&\xi^6\sum_{\gamma_{r_0(^\star x_{i,j}), (x_{k,l}^\star)r_0r_3, h}\ne 0}\gamma_{r_0(^\star x_{i,j}), (x_{k,l}^\star)r_0r_3, h}C_{h}C_{r_1}+\triangle+\Box,\\
=&\xi^6\sum_{\gamma_{r_0(^\star x_{i,j}), (x_{k,l}^\star)r_0r_3, r_0(^\star x_{m,n}^\star)r_0r_3}\ne 0}\gamma_{r_0(^\star x_{i,j}), (x_{k,l}^\star)r_0r_3, r_0(^\star x_{m,n}^\star)r_0r_3}C_{r_0(^\star x_{m,n}^\star)r_0r_3}C_{r_1}+\triangle+\Box,\\
=&\xi^6\sum_{\gamma_{r_0(^\star x_{i,j}), (x_{k,l}^\star)r_0r_3, r_0(^\star x_{m,n}^\star)r_0r_3}\ne 0}\gamma_{r_0(^\star x_{i,j}), (x_{k,l}^\star)r_0r_3, r_0(^\star x_{m,n}^\star)r_0r_3}C_{r_0(^\star x_{m,n}^\star)r_0r_3r_1}+\triangle+\Box\\
=&\xi^6\sum_{\gamma_{r_0(^\star x_{i,j}), (x_{k,l}^\star)r_0r_3, r_0(^\star z^\star)r_0r_3}\ne 0}\gamma_{r_0(^\star x_{i,j}), (x_{k,l}^\star)r_0r_3, r_0(^\star z^\star)r_0r_3}C_{r_0(^\star z^\star)r_0r_3r_1}+\triangle+\Box,
\end{aligned}\end{equation}
where $\triangle\in H$ is a sum of elements in $\Gamma'_{013}\cap\Gamma_{02}^{-1}$ whose coefficient is of degree less than 6.

So we have shown (42) and (35) is valid afterwards.

For (36), by 1.4(d), we have 
$\gamma_{(x_{i,j}^\star)r_0, r_0(^\star x_{k,l}), h_5}=\gamma_{r_0(^\star x_{k,l}),  h_5^{-1}, r_0(^\star x_{i,j})}$, and $h_5\in\Gamma_{012}\cap\Gamma_{012}^{-1}$.

By (18) and 1.4(b)(g), we have 
$\gamma_{r_0(^\star x_{k,l}),  h_5^{-1}, r_0(^\star x_{i,j})}=\gamma_{^\star x_{k,l},  h_5^{-1}, ^\star x_{i,j}}=\gamma_{x_{k,l},  h_5^{-1}, x_{i,j}}=\gamma_{ x_{i,j}, x_{k,l},  h_5}.$

For (37), it is deduced by (18), (36) and 1.4(b)(g).
For (38), it is deduced by (37) and 1.5(b).
For (39), it is deduced by (28), (37) and 1.4(b).

\medskip

To prove assertion (c) and (d) of Proposition 3.9A, by (a) and (b) of Proposition 3.9A we need to check the following identity:
for any $r_3x\in\Gamma_{012}\cap\Gamma_{03}^{-1}, y, r_3z\in W$,
 we have $x\in\Gamma_{012}\cap\Gamma_{02}^{-1}$ and
\begin{align}
\gamma_{r_3x, y, r_3z}=\gamma_{x, y, z}.\end{align}

Note by 1.4(c) we have $\gamma_{x, y, z}\ne 0$ only if $z\in\Gamma_{02}^{-1}$, so $z$ is the first element of the left string with respect to $\{r_2,r_3\}$ and $r_3x$ is a second element of the left string with respect to $\{r_2,r_3\}$.

By 1.5(a), we get $\gamma_{r_3x, y, r_3z}=\gamma_{x, y, z}+\gamma_{x, y, r_2r_3z}$. Since $r_2r_3z\in(\Gamma'_{02})^{-1}$, by 1.4(c), we get $\gamma_{x, y, r_2r_3z}=0$ and (44) is true.

\medskip

To prove assertion (e) and (f) of Proposition 3.9A, by (c) and (d) of Proposition 3.9A we need to check the following identity:
for any $x\in\Gamma_{012}\cap\Gamma_{03}^{-1}, y, z\in W$, we have
\begin{align}
\gamma_{r_1x, y, r_1z}=\gamma_{x, y, z}.\end{align}

Note by 1.4(c) we have $\gamma_{x, y, z}\ne 0$ only if $z\in\Gamma_{03}^{-1}$.

By (27),1.4(b) and (c) (d) of Proposition 3.9A , we get

$\begin{aligned}
C_{r_1x}C_y
=&C_{r_1}C_xC_y+\Box\\
=&\xi^6\sum_{\gamma_{x, y, z}\ne 0}\gamma_{x, y, z}C_{r_1}C_z+\triangle+\Box\\
=&\xi^6\sum_{\gamma_{x, y, z}\ne 0}\gamma_{x, y, z}C_{r_1z}+\triangle+\Box,\end{aligned}$

where $\triangle\in H$ is a sum of elements in $\Gamma_{012}\cap(\Gamma'_{013})^{-1}$ whose coefficient is of degree less than 6.

So (45) is valid.

\medskip

To prove assertion (g) and (h) of Proposition 3.9A, by Proposition 3.6A we need to check the following identity:
for any $r_3x\in\Gamma_{02}\cap\Gamma_{03}^{-1}, y, z\in W$ we have
$x\in\Gamma_{02}\cap\Gamma_{02}^{-1}$,and
\begin{align}
\gamma_{r_3x, y, r_3z}=\gamma_{x, y, z}.\end{align}

Note $\gamma_{x, y, z}\ne 0$ forces $y, z\in\Gamma_{02}^{-1}$.

By (23), Proposition 3.6A and 1.4(b), we get

$\begin{aligned}
C_{r_3x}C_y
=&C_{r_3}C_xC_y+\Box\\
=&\xi^6\sum_{\gamma_{x, y, z}\ne 0}\gamma_{x, y, z}C_{r_3}C_z+\triangle+\Box\\
=&\xi^6\sum_{\gamma_{x, y, z}\ne 0}\gamma_{x, y, z}C_{r_3z}+\triangle+\Box,\end{aligned}$

where $\triangle\in H$ is a sum of elements in $\Gamma_{03}^{-1}$ whose coefficient is of degree less than 6.

So (46) is valid. 

\medskip

To prove assertion (i) and (j) of Proposition 3.9A, by (g) and (h) of  Proposition 3.9A we need to check the following identity:
for any $x\in\Gamma_{02}\cap\Gamma_{03}^{-1}, y, z\in W$, we have
\begin{align}
\gamma_{r_1x, y, r_1z}=\gamma_{x, y, z}.\end{align}

Note $\gamma_{x, y, z}\ne 0$ forces $y\in\Gamma_{02}^{-1}, z\in\Gamma_{03}^{-1}$.

By (27), (g) and (h) of  Proposition 3.9A and 1.4(b), we get

$\begin{aligned}
C_{r_1x}C_y
=&C_{r_1}C_xC_y+\Box\\
=&\xi^6\sum_{\gamma_{x, y, z}\ne 0}\gamma_{x, y, z}C_{r_1}C_z+\triangle+\Box\\
=&\xi^6\sum_{\gamma_{x, y, z}\ne 0}\gamma_{x, y, z}C_{r_1z}+\triangle+\Box,\end{aligned}$

where $\triangle\in H$ is a sum of elements in $(\Gamma'_{013})^{-1}$ whose coefficient is of degree less than 6.

So (47) is valid.

\medskip

To prove assertion (i) and (j) of Proposition 3.9A, by Proposition 3.7A we need to check the following identity:
for any $x\in\Gamma_{03}\cap\Gamma_{03}^{-1}, y, z\in W$, we have
\begin{align}
\gamma_{r_1x, y, r_1z}=\gamma_{x, y, z}.\end{align}

Note $\gamma_{x, y, z}\ne 0$ forces $y, z\in\Gamma_{03}^{-1}$.

By (27), Proposition 3.7A and 1.4(b), we get

$\begin{aligned}
C_{r_1x}C_y
=&C_{r_1}C_xC_y+\Box\\
=&\xi^6\sum_{\gamma_{x, y, z}\ne 0}\gamma_{x, y, z}C_{r_1}C_z+\triangle+\Box\\
=&\xi^6\sum_{\gamma_{x, y, z}\ne 0}\gamma_{x, y, z}C_{r_1z}+\triangle+\Box,\end{aligned}$

where $\triangle\in H$ is a sum of elements in $(\Gamma'_{013})^{-1}$ whose coefficient is of degree less than 6.

So (48) is valid.

The proof is completed.
\end{proof}

\medskip

 {\bf Corollary 3.9B.} For any left cell $\Gamma\in Y_i,\ \Theta\in Y_j$ with $i\ne j$, we have a bijection
 $$\pi_{\Gamma\cap\Theta^{-1}}: \Gamma \cap\Theta^{-1}\to\text{the set of isomorphism classes of irreducible $F_c$-v.b. on}\ \Theta\times \Gamma,$$
 inducing an isomorphism of $\mathbb Z$-modules
  $$\pi_{\Gamma \cap\Theta^{-1}}: J_{\Gamma \cap\Theta^{-1}}\to K_{F_c}(\Theta\times \Gamma)\cong\text{Rep}\ (Sp_4(\mathbb C)\times\mathbb Z/2\mathbb Z),$$
  such that take any $x\in\Gamma \cap\Theta^{-1}$, $\Phi\in Y_1\cup Y_2\cup Y_3\cup Y_4$, we have

(a) $\pi_{\Phi\cap\Theta^{-1}}(t_xt_y)=\pi_{\Gamma\cap\Theta^{-1}}(t_x)\pi_{\Phi\cap\Gamma^{-1}}(t_y), \text{for\ any}\ y\in \Phi\cap\Gamma^{-1}$,

(b) $\pi_{\Gamma\cap\Phi^{-1}}(t_zt_x)=\pi_{\Theta\cap\Phi^{-1}}(t_z)\pi_{\Gamma\cap\Theta^{-1}}(t_x), \text{for\ any}\ z\in \Theta\cap\Phi^{-1}$.

  The inverse gives a bijection
  $$\pi_{\Theta\cap\Gamma ^{-1}}: \Theta\cap\Gamma ^{-1}\to\text{the set of isomorphism classes of irreducible $F_c$-vector bundles on}\ \Gamma\times \Theta,$$
  inducing an isomorphism of $\mathbb Z$-modules
 $$\pi_{\Theta\cap\Gamma ^{-1}}: J_{\Theta\cap\Gamma ^{-1}}\to K_{F_c}(\Gamma\times \Theta)\cong\text{Rep}\ (Sp_4(\mathbb C)\times\mathbb Z/2\mathbb Z),$$
   such that take any $x\in\Theta\cap\Gamma ^{-1}$, $\Phi\in Y_1\cup Y_2\cup Y_3\cup Y_4$, we have

(a) $\pi_{\Phi\cap\Gamma^{-1}}(t_xt_y)=\pi_{\Theta\cap\Gamma ^{-1}}(t_x)\pi_{\Phi\cap\Theta^{-1}}(t_y), \text{for\ any}\ y\in \Phi\cap\Theta^{-1}$,

(b) $\pi_{\Theta\cap\Phi^{-1}}(t_zt_x)=\pi_{\Gamma\cap\Phi^{-1}}(t_z)\pi_{\Theta\cap\Gamma ^{-1}}(t_x), \text{for\ any}\ z\in \Gamma\cap\Phi^{-1}$.
\begin{proof}
It follows Lemma 3.2 and subsection 3.3.
  \end{proof}

Combining Proposition 3.9A and Corollary 3.9B, we have proved Theorem 3.4B.

\bigskip

By Theorem 3.4A and Theorem 3.4B, the proof of Theorem 3.4 now is completed.

\bigskip

\noindent{\bf Acknowledgement:} Part of the work was done during
my visit to the Academy of Mathematics and Systems Science, Chinese Academy of Sciences. I am very grateful to the AMSS for hospitality and financial supports. And I would like to thank Nanhua Xi for useful discussions.

\end{spacing}
\end{document}